# A New Slant on Lebesgue's Universal Covering Problem

By Philip Gibbs

Abstract

Lebesgue's universal covering problem is re-examined using computational methods. This leads to conjectures about the nature of the solution which if correct could provide a blueprint for a complete solution. Empirical lower bounds for the minimal area are computed using different hypotheses based on the conjectures. A new upper bound of 0.844113 for the area of the minimal cover is derived improving previous results. This method for determining the bound is suggested by the conjectures and computational observations but is proved independently of them. The key innovation is to modify previous best results by removing corners from a regular hexagon at a small slant angle to the edges of the dodecahedron used before. Simulations indicate that the minimum area for a convex universal cover is likely to be around 0.84408.

## 1. Introduction

It is 100 years since Henri Lebesgue posed a curious problem about universal covers in the plane for sets of diameter one. The diameter of a bounded set of points $A$ is the supremum of the distance $d(x, y)$ for $x, y \in A$. In the original problem written in a letter to Pál in 1914, Lebesgue asked for the convex set of minimum area that contains a subset congruent to any set of diameter one [1]. By the definition of congruence this means that the shape $A$ can be translated, rotated or reflected to fit in the cover. This is known as **Lebesgue's Universal Covering Problem** or the **Lebsegue Minimal Problem**

A variation is to look for the minimal area of non-convex universal covers. For the convex case there must be a unique universal cover that achieves the minimum area by the Blaschke selection theorem [2]. It is not known if a minimal universal cover exists for shapes that may not be convex but the infinum of the areas of all covers is well defined.

In 1920 Pál showed that a regular hexagon circumscribing a circle of unit diameter is a universal cover. This set an upper limit for Lebesgue's minimal area at $x \leq \frac{\sqrt{3}}{2} = 0.86602540$ He also showed that two corners can be removed from the hexagon down to the sides of a regular dodecagon inside the hexagon to make a smaller universal cover. This set a better upper limit $x \leq 2 - \frac{2}{\sqrt{3}} = 0.84529946$ [3]. Then in 1936 Sprague showed that a smaller piece bounded by two circular arcs could be removed from another corner [4] reducing the upper bound to

$$x \leq 0.844137708435197570894066994 \quad (1)$$

Finally Hansen between 1975 and 1992 [5-7] removed two further small pieces from this shape. This reduced the upper bound on the area by only 1.8738 x $10^{-11}$ to

$$x \leq 0.844137708416458789749982 2 \quad (2)$$

Hansen's argument made use of the freedom to reflect shapes which had not been required by Pál or Sprague and produced a covering without bilateral symmetry.

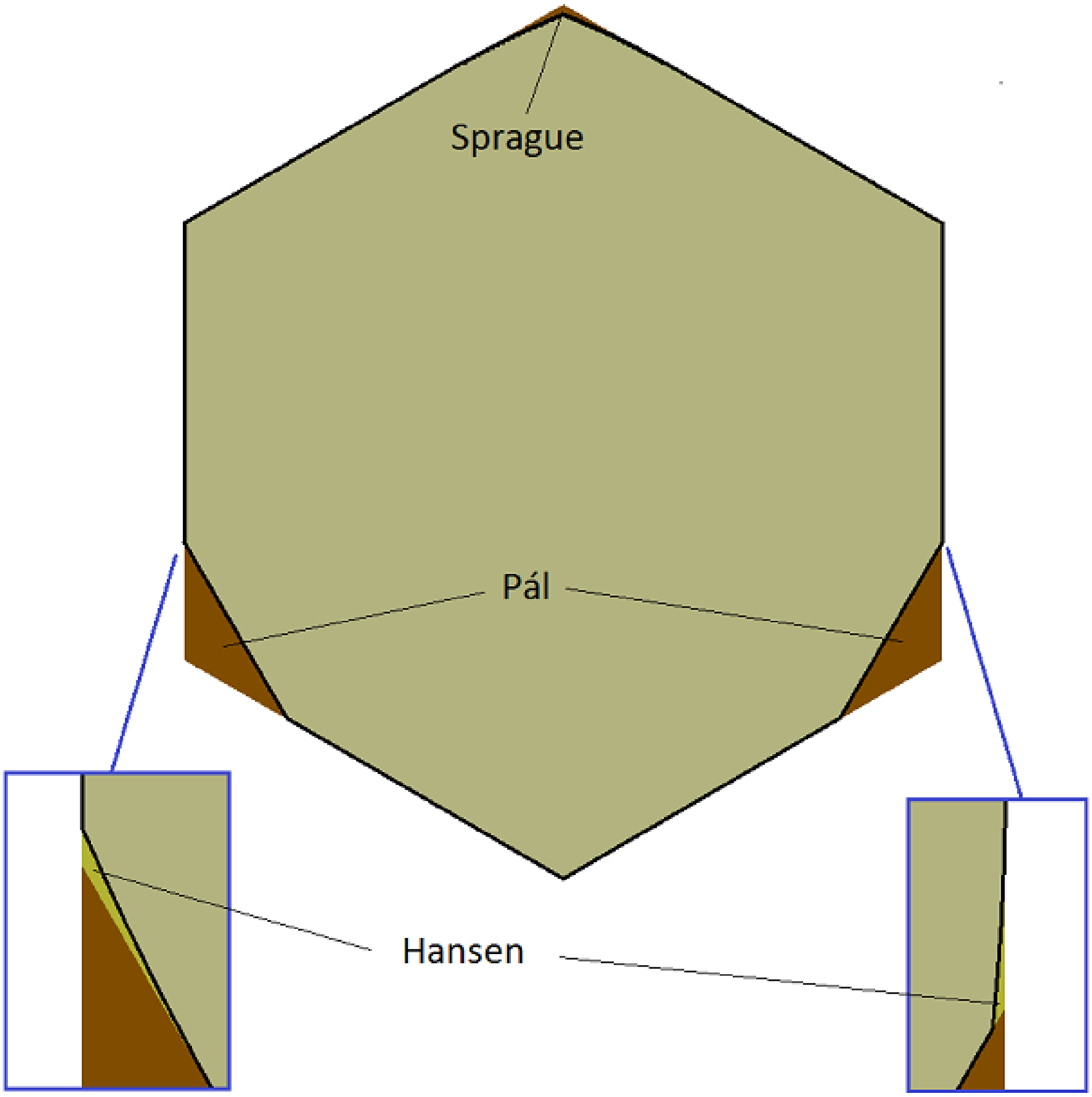

Figure 1: minimal convex universal cover from previous results

Another small reduction to Sprague's covering was made in 1980 by Duff using a non-convex cover giving an area of 0.84413570 for the non-convex variation on the problem. This reduction also made use of reflections [8]. Duff considered another variation of the problem proposed by Rennie in 1977 [9] in which only translations of the shapes are permitted [10]. This is useful as a simpler but still unresolved variant whose solution would help understand the harder problem set by Lebesgue.

Proofs for lower bounds have so far been much further from the likely minimal answer but nevertheless they are worthy results. Elekes was able to compute the minimal area for a circle and all Reuleaux polygons for which the number of sides is a power of three [11]. The best lower bound by Brass and Sharifi sets a lower bound of 0.832 by combining a circle, triangle and pentagon of diameter one [12].

More related problems are discussed in the book by Brass, Moser and Pach [13]

The present work will describe new computational methods which provide better lower bounds under some assumptions of the shape containing the cover. This study has also led to a new rigorous argument improving the upper bound of Hansen for convex cover.

## 2. Simulated Annealing

A basic and easily implemented computational technique used for minimisation problems with many variables is simulated annealing. With a finite number of shapes it is possible to run a simulation in which the shapes are allowed to translate and rotate randomly as if in thermal motion at a given temperature. The combined area of the shapes is used as the energy function.

The simulated temperature can be dropped slowly so that the shapes freeze into a local minimum of the area. If the cooling is done slowly and the simulation is repeated many times with different randomisations, it is possible to empirically locate the minimum area. In Lebesgue's problem it is also allowed to reflect shapes to form the minimum cover and the simplest method to accommodate this is to simulate $2^n$ different combinations of reflections where $n$ is the number of shapes in the simulation without bilateral symmetry.

Thermal motion at a given temperature can be simulated with the correct thermal distribution using a Metropolis-Hastings algorithm. However, this would require the size of the movements to be matched in a carefully controlled way to the temperature as it decreased; otherwise the efficiency of the method is lost. A simpler method is to use a biased random walk in which changes that decrease the area are accepted, while those that increase the area are kept with a fixed probability (e.g. $p$=0.3). With this adaptation the effective temperature decreases naturally as the size of trial changes decreases while the selection criterion remains efficient.

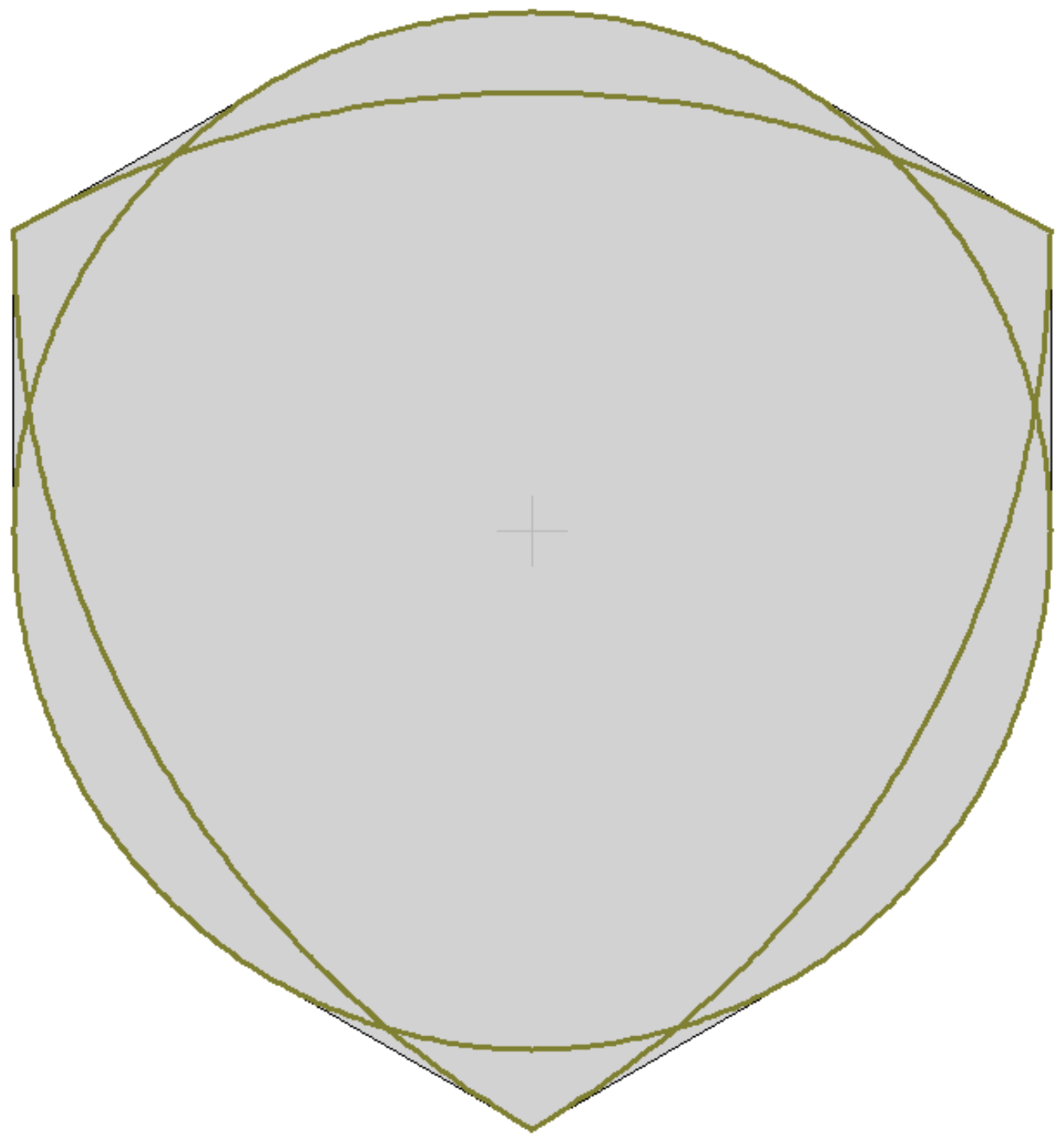

**Figure 2: Circle and Reuleaux Triangle**

The optimal shapes to use are curves of constant width. Indeed any shape of diameter one is contained within a curve of constant width equal to one [14], so it is sufficient to consider only curves of constant width in the covering problem. Reuleaux polygons are the most easily constructed curves of constant width and empirically they appear to be most effective in maximising the minimum area of a cover for a given number of shapes. The curve of constant width with the maximum area is the circle, while the Reuleaux triangle has minimum area (figure 2). Together these two shapes maximise the area of the minimum convex cover for two shapes.

In general a Reuleaux polygon is formed from an odd number *n* of circular arcs of radius one whose centres are placed on a star polygon whose sides all have length one. It is sufficient to specify *n*-3 consecutive angles between these sides which form the diameters of the Reuleaux shape. The remaining point and three angles can be determined by triangulation from the two loose ends. All the angles must be less than 60 degrees and will add up to 180 degrees. The necessary requirement that all diagonals are less than or equal to one poses further constraints on the allowed range of angles.

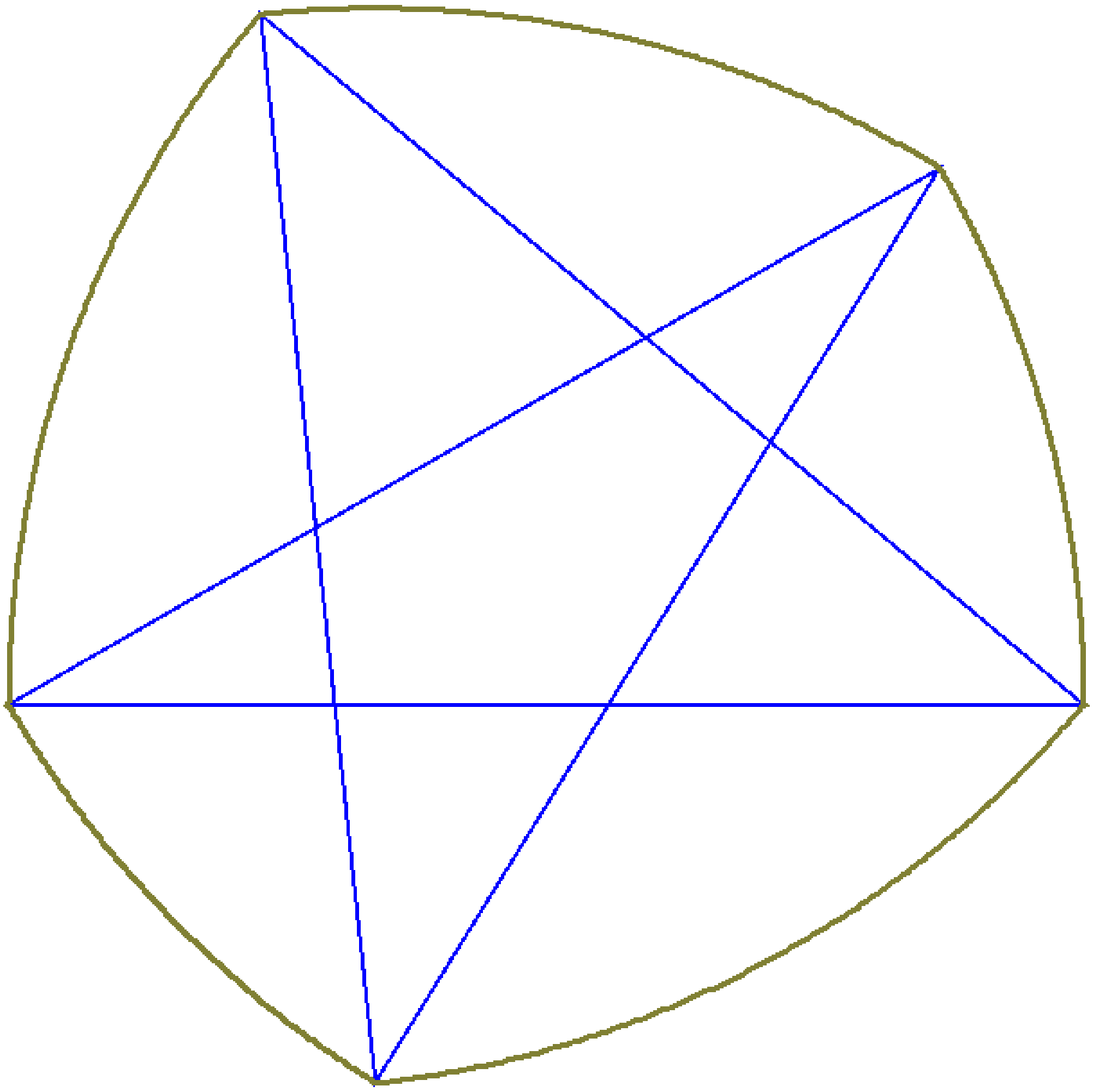

Figure 3: construction of a Reuleaux pentagon

In the simulations a small number of Reuleaux polygons were used. For simplicity the arcs were replaced by polygons with small sides ensuring that the diameter of the shape remained equal to one. The area of the convex hull of the shapes was calculated using a simple wrapping algorithm.

More efficient algorithms exist (see http://en.wikipedia.org/wiki/Convex_hull_algorithms) and it would be possible to implement them for arc-polygons rather than straight edge polygons. In this way the simulation could be made much more efficient. The simulation was implemented as a Java applet.

As an example of the results this shows the final state of a simulation with five shapes: a circle, and regular Reuleaux polygons with 3,5,7,and 9 sides.

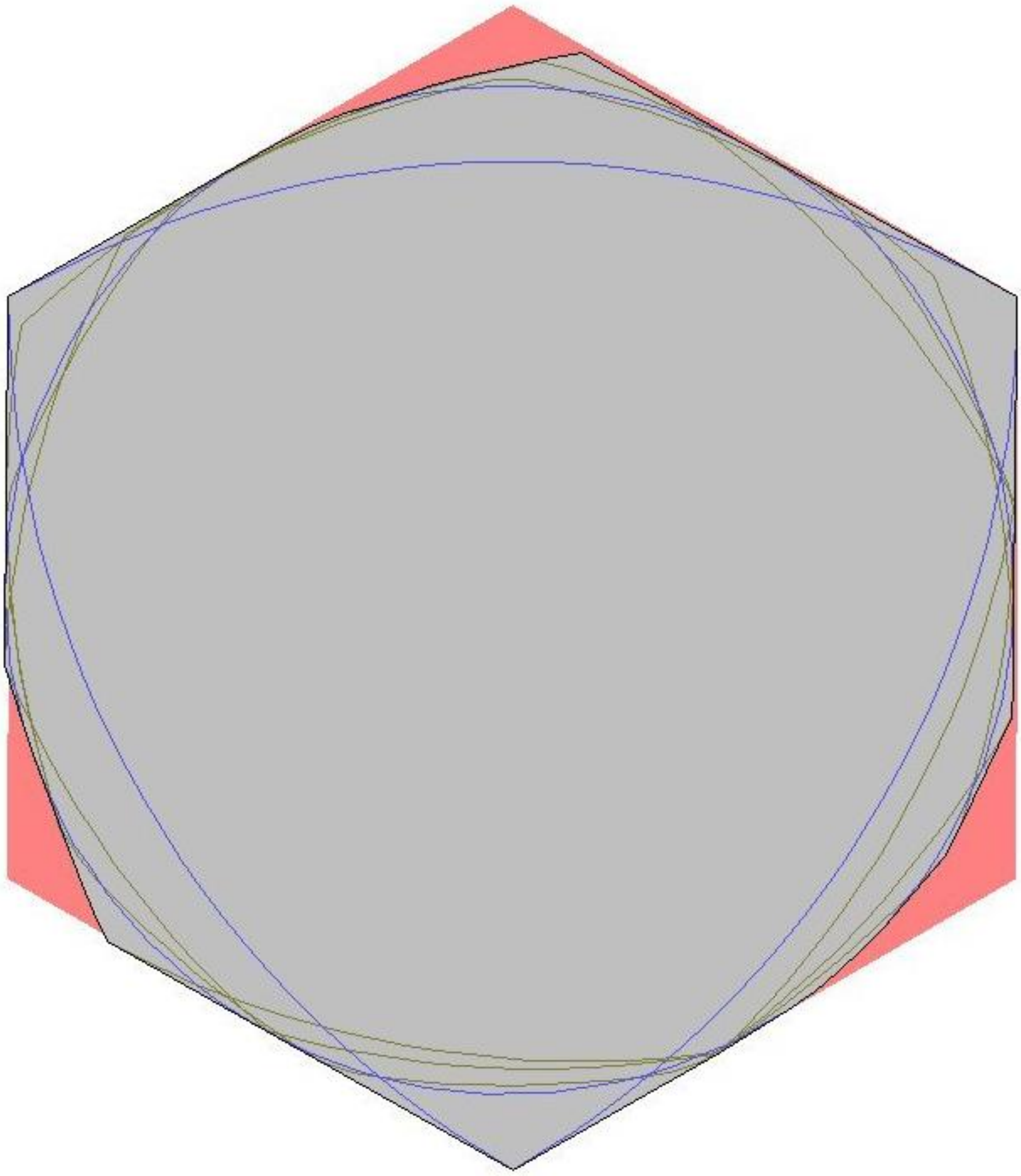

**Figure 4: results of a simulated annealing for 5 shapes**

The minimum area for any given set of shapes of diameter one is a lower bound on the minimal convex cover. The answer obtained by simulated annealing can only be regarded as an upper bound on such a lower bound, but with multiple simulations and a small number of shapes the best area can be found with some certainty and accuracy and can be regarded as an empirical lower bound. For these shapes the minimal convex area found was 0.83699098.

As the number of shapes increases, the number of local minima makes it hard to find the true minimum of area. However, a useful observation seen in all cases tried is that the minimal cover appears to sit within an irregular hexagon whose opposite sides are parallel and distance one apart. In the figure above a regular hexagon has been drawn as a background reference shape making it easier to see that the edges of the hexagon that contain the shapes are at a slight angle to the regular shape.

This can be expressed as a conjecture:

*The minimum area of a convex cover for any set of shapes of constant width one can be contained in a hexagon whose three pairs of opposite sides are parallel and a distance one apart*

If the set of shapes includes a circle then the sides of the hexagon would be tangent to the circle, the opposing pairs would be equal in length and the hexagon would have a 180 degree rotational symmetry

## 3. The Hexagon Hypothesis

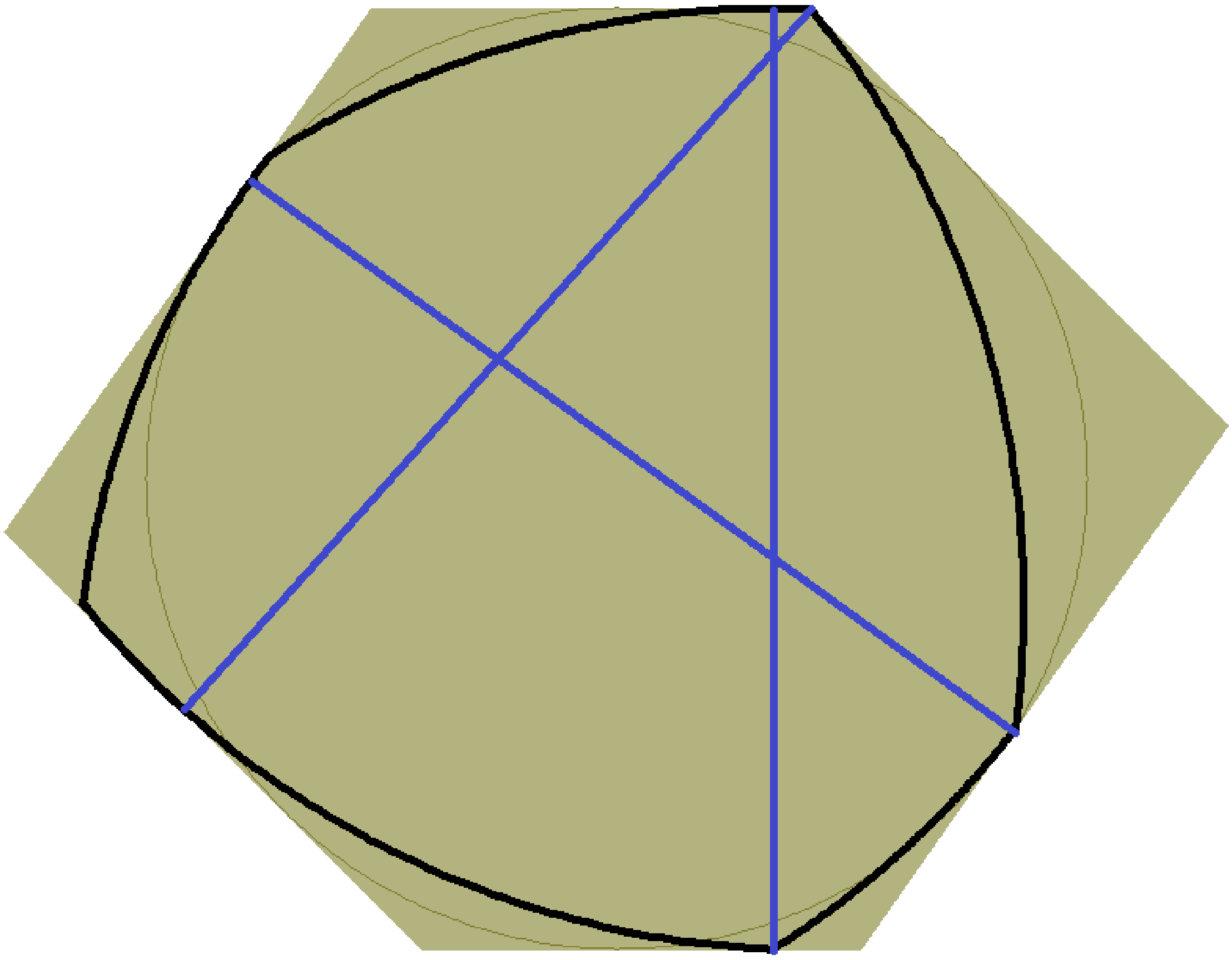

Figure 5: curve of constant width in a hexagon

Pál showed that a regular hexagon circumscribed around a circle of unit diameter is a universal cover, but this can be generalised to any hexagon circumscribed around the circle such that opposite sides are parallel. For any such circumscribed parallel hexagon there is a minimal area for covers that are bounded by that hexagon. The regular hexagon has the least total area but that does not necessarily mean that the minimum area for a cover contained in a circumscribed parallel hexagon is contained in the regular hexagon.

The results from simulated annealing suggest that a cover of minimal area for a set of shapes that include the circle can be contained in such a hexagon but in absence of a proof we cannot be certain. It may be that the covers generated by annealing are only approximately contained in the hexagon or that with more shapes they would take a different form. Nevertheless we can work on the basis of a "hexagon hypothesis" and seek the minimum area that is contained in such a hexagon. This would be at least a good upper bound for the true answer.

First we prove that a circumscribed parallel hexagon is always a universal cover. Define the support function $s(\theta)$ for any shape relative to a point C as a function of an angle $\theta$ as the distance from C to a tangent to the shape normal to the direction $\theta$. Without using the notion of tangent the support function is defined for any compact set $S$ in the plane as

$$s(\theta) = \sup((P - C).\boldsymbol{n} | P \in S \tag{3}$$

Where $\boldsymbol{n}$ is the unit vector in direction $\theta$ i.e. $\boldsymbol{n} = (\cos\theta, \sin\theta)$

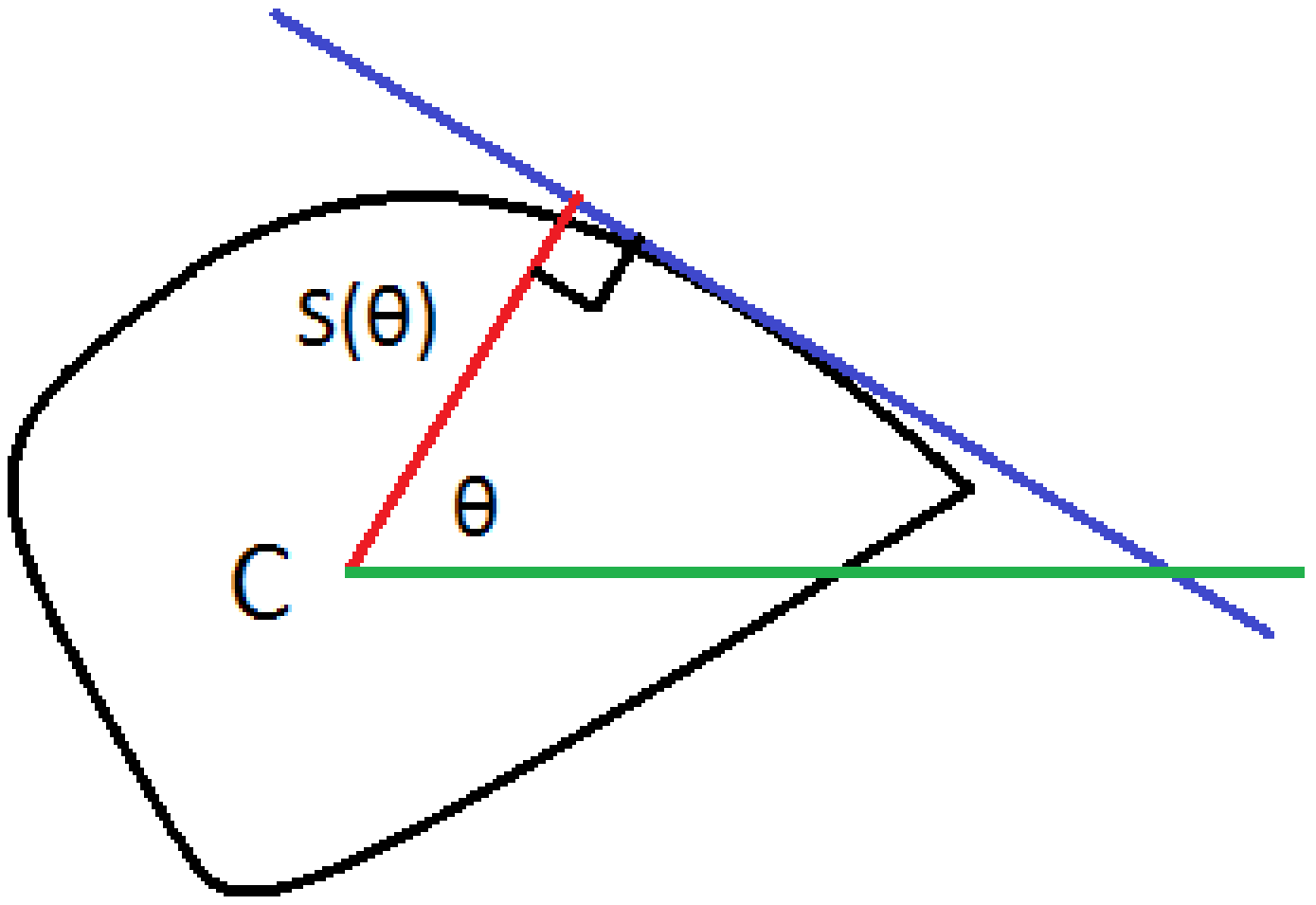

Figure 6: support function for a closed curve

The condition for the shape to be a curve of constant width one is that $s(\theta) + s(\theta + \pi) = 1$. This condition is independent of the choice of C. An offset function defined by

$$o(\theta) = (s(\theta) - s(\theta + \pi))/2 \tag{4}$$

Gives the offset of the midpoint of the diameter in direction $\theta$ from C in the same direction. For a curve of constant width it satisfies $o(\theta + \pi) = -o(\theta)$. Now consider the offsets $o_1, o_2, o_3$ in three directions given by the unit normals $\boldsymbol{n_1}, \boldsymbol{n_2}, \boldsymbol{n_3}$ to the sides of a circumscribed parallel hexagon centred at a point H. If H is allowed to move so that the vector from C to H is $\boldsymbol{x}$ then the shape would be covered by the hexagon when $\boldsymbol{x}.\boldsymbol{n_i} = o_i,\ i = 1, 2, 3$. With the ansatz $\boldsymbol{x} = a\widehat{\boldsymbol{n}}_{\boldsymbol{1}} + b\widehat{\boldsymbol{n}}_2$ where $\widehat{\boldsymbol{n}}_i$ is a unit vector at 90 degrees to $\boldsymbol{n}_i$ , the first two equations are solved by

$$a = \frac{x.n_2}{n_1 \wedge n_2}, b = \frac{x.n_1}{n_2 \wedge n_1} \tag{5}$$

The third equation is equivalent to

$$t = o_1(\boldsymbol{n}_3 \wedge \boldsymbol{n}_2) + o_2(\boldsymbol{n}_1 \wedge \boldsymbol{n}_3) + o_3(\boldsymbol{n}_2 \wedge \boldsymbol{n}_3) = 0 \tag{6}$$

This equation is independent of C since the requirement that the hexagon can be translated to cover the shape is independent of C. To attempt to solve this condition and fit the curve inside the hexagon we rotate the hexagon by an angle $\theta$ so that $t$ is a function $t(\theta)$. The angles at which the shape fits the hexagon are the roots of $t(\theta) = 0$. That such roots exist follows from the identity $t(\theta + \pi) = -t(\theta)$ and the fact that $t(\theta)$ is continuous. This proves that any curve of constant width can be made to fit inside any circumscribed parallel hexagon and therefore that such a shape is a universal cover.

Furthermore the support function for a Reuleaux polygon is easily calculated so a function to find all the roots of $t(\theta)$ can be programmed using recursive bisection making use of the fact that the derivative $t'(\theta)$ is bounded to eliminate ranges that do not contain roots. For computational purposes the implementation is robust, efficient and accurate. This offers the possibility of finding the minimum area for a cover of a given set of shapes inside a given hexagon. For each shape there are a finite number of ways it can be placed inside the hexagon allowing translations, rotations and reflections and these can be computed. It is then sufficient to search through all the ways that all the shapes can be placed inside the hexagon to find the minimum area that covers them, either using a convex cover or the a non-convex cover to compute the area.

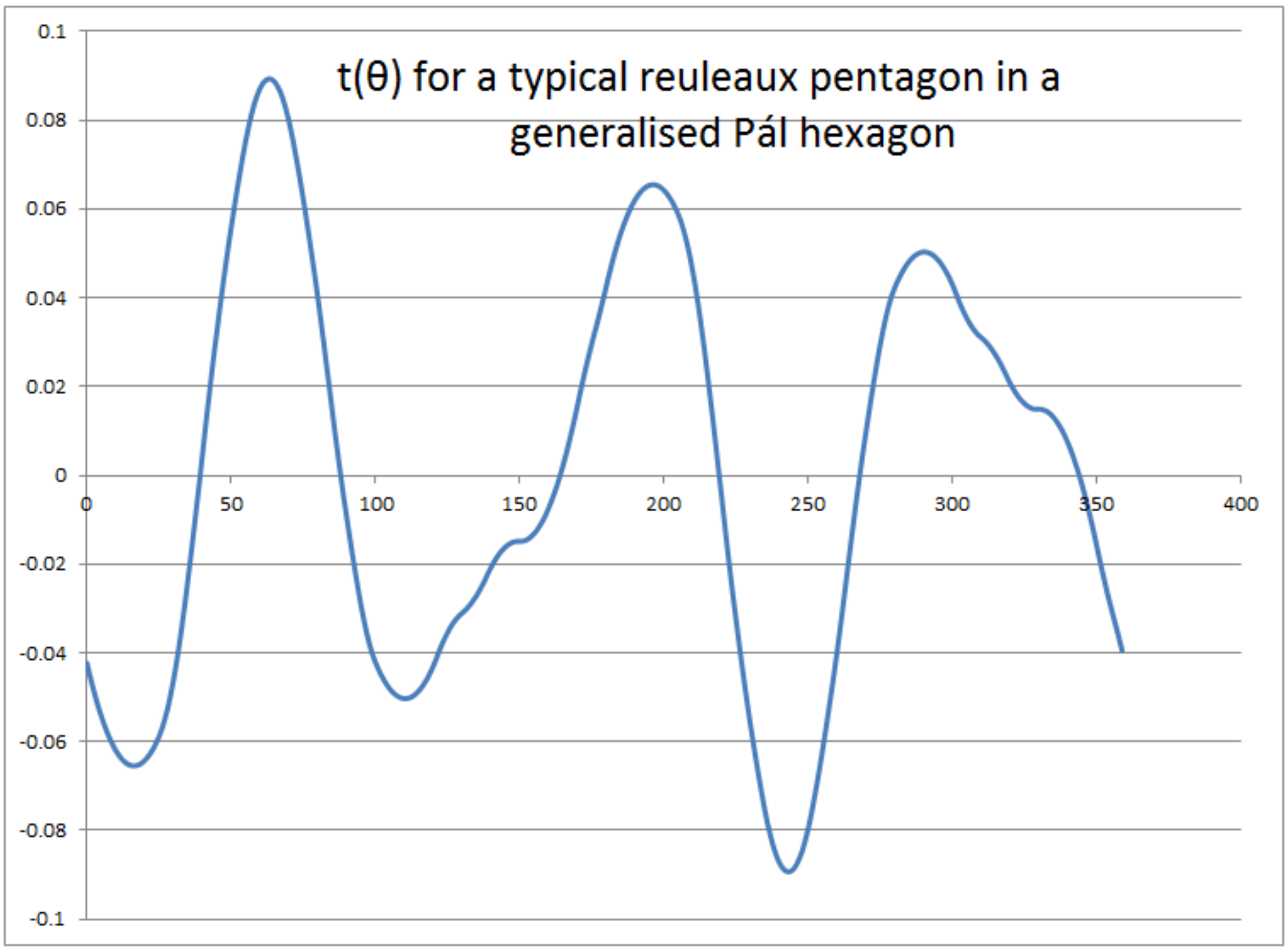

Figure 7 offset function for a typical Reuleaux pentagon

For a typical curve of constant width we find that $t(\theta)$ has six roots, though it can have fewer or an unlimited number more, the chart below shows $t(\theta)$ for three symmetrical pentagons where two anlges are set to 35, 36 and 37 degrees. The numbers of roots are 18, 30 and 6 respectively.

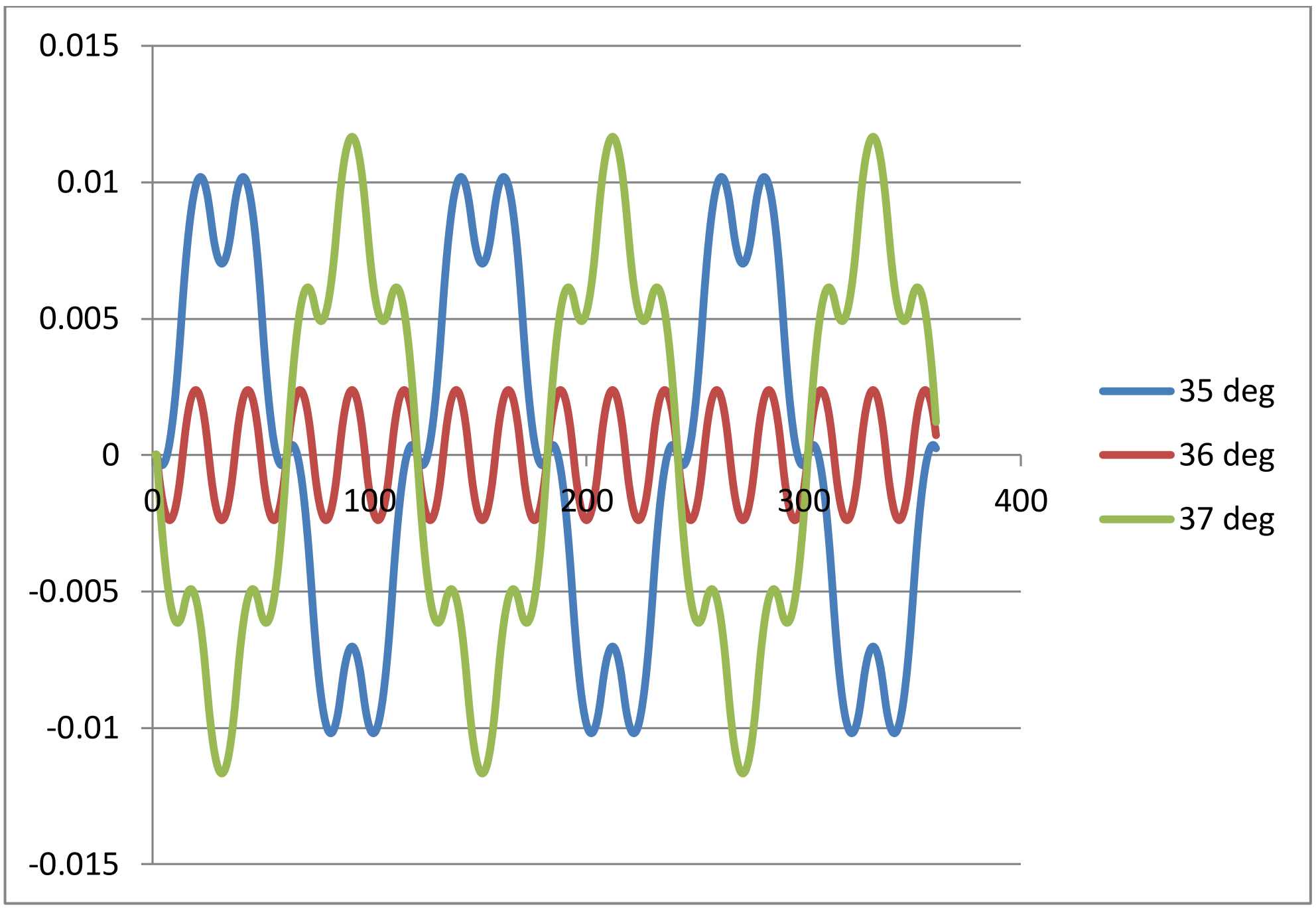

Figure 8: offset funtions for three Reuleaux pentagons

Together with the reflections if it is asymmetrical, this means that there are typically $12^N$ areas to calculate to find the minimum area. This would limit the search to just a few shapes were it not for the possibility of pruning the search tree using the fact that the minimum area can only grow as more shapes are added. At any point in the search the minimum area known so far can be used as a cutoff to prune further addition of shapes. In practice it pays to calculate all possible areas for the addition of any one further shape at any point in the search, and choose to add the shape that increases its minimum additional area by more than any other shape. This takes us to the optimum minimum more quickly and means that we reach a point where any branch of the tree can be pruned more quickly. Using these optimisations it has been possible to find the minimum area for well over a hundred shapes inside a given circumscribed parallel hexagon.

In an implementation of this algorithm shapes are selected to be added one at a time from a large pool of Reuleaux polygons. Each time the shape which appears to be least well covered by the minimum area defined by previous shapes is added, then the new least cover is computed. This is repeated until the area appears to converge. At any time this gives a rigorous lower bound for the minimum cover included in a given hexagon. If it exceeds the known upper bound for the minimal universal cover then that shape can be ruled out as a container for a minimum cover

A circumscribed parallel hexagon can be specified by giving the angles between its edges. Only two angles are required since opposite sides are parallel and the three remaining different angles add up to 180 degrees. Using the search algorithm we can explore the two dimensional space of circumscribed parallel hexagons parameterised by these two angles. A graph of the lower bound on the area using up to 40 shapes in the special case where the two angles are equal can be plotted as follows (The regular hexagon is the case where the angles are 60 degrees) The horizontal line is the known upper limit so we can see that there is only a possibility of finding a lower area when the angles are within 0.5 degree of the regular hexagon.

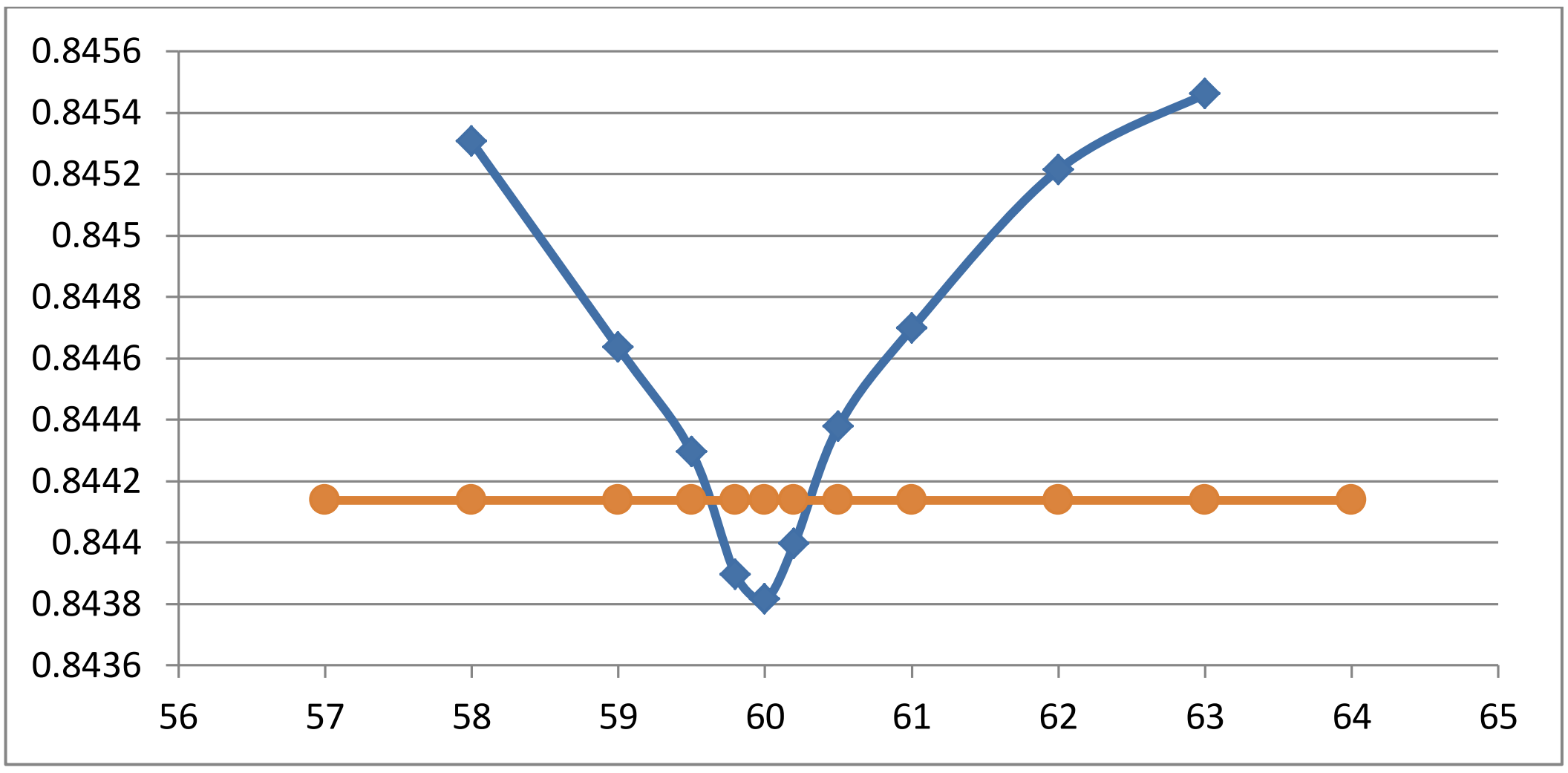

Figure 9: lower bound for symmetric hexagon

A plot closer to this point using up to 60 shapes shows that the minimum appears to be at or very close to the case of the regular hexagon. With 80 shapes the lower bound for the regular hexagon is 0.843961

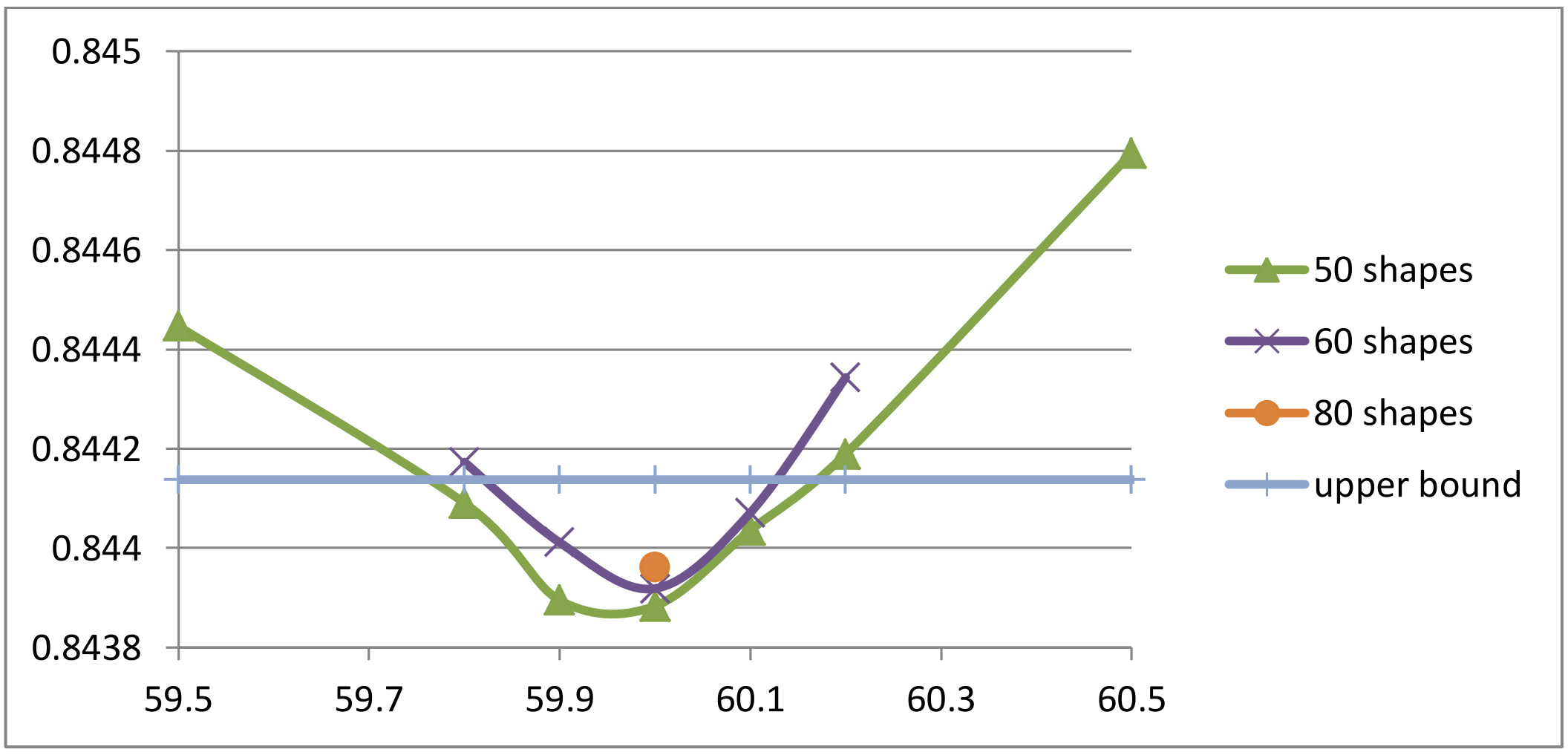

Figure 10: lower bound for symmetric hexagon

More of the space can be covered e.g. by taking one angle to be 60 degrees and varying a second resulting in this plot

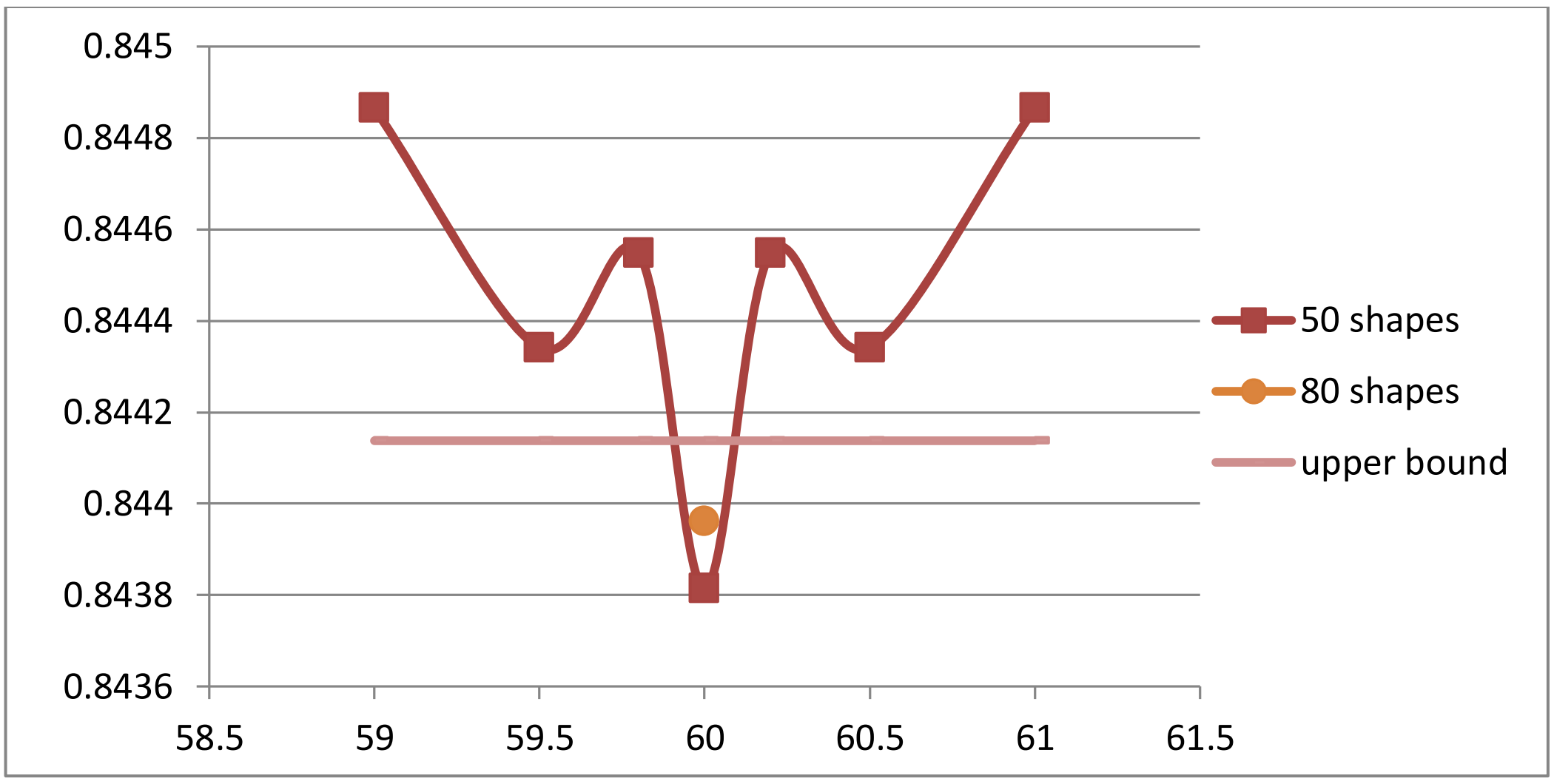

Figure 11: lower bound for asymmetric hexagon with one angle of 60 degrees

Once again the conclusion is that the regular hexagon is the best case to investigate further although it is possible that the minimal cover is inside a hexagon that is irregular, but differing only slightly from regular.

The conclusion a second conjecture is proposed

*The minimal universal cover within a hexagon with the property that opposite sides are parallel and distance one apart is achieved in the case of the regular hexagon*

## 4. Modified Pál Hypothesis

For the regular case the minimal cover with 80 shapes looks like this

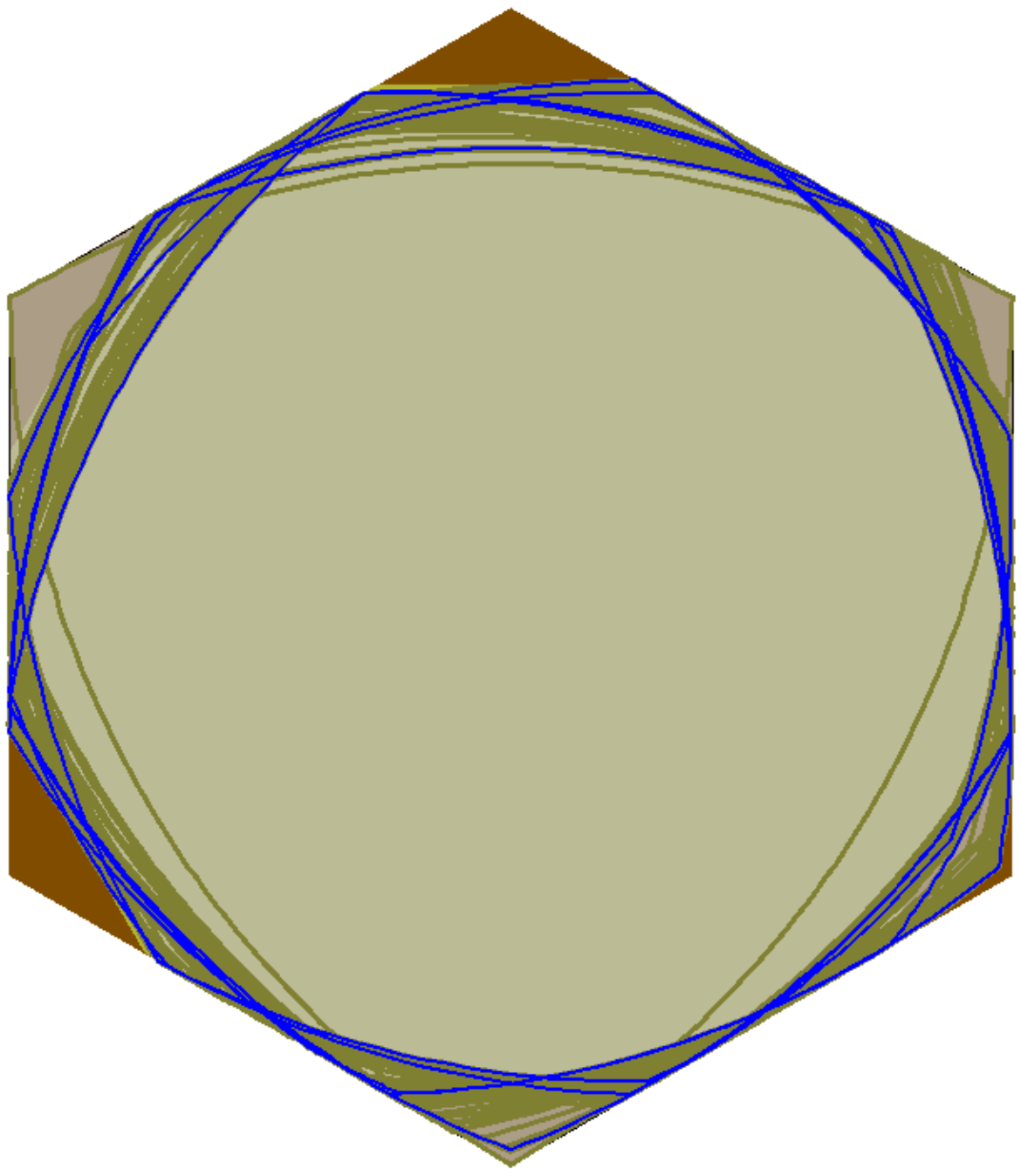

Figure 12: minimal cover inside regular hexagon for 80 shapes

If we zoom in on one of the cut-off corners we find that the minimal cover differs from Pál's solution in that the edge of the corner cut off is at a small angle to the edge of the regular dodecahedron.

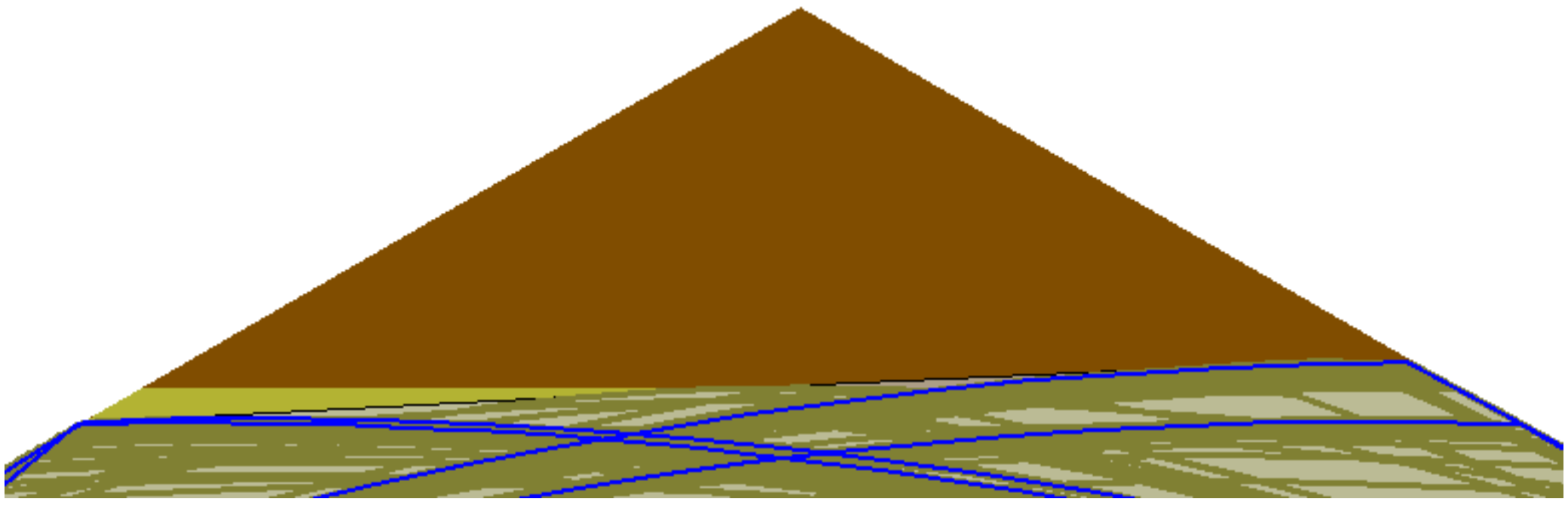

Figure 13: minimal cover in hexagon showing slanting cutoff

Pál demonstrated that two corners of a regular hexagon circumscribing a circle can be removed to construct a smaller universal cover. The cuts were tangent to the unit circle forming two sides of a regular dodecagon inside the hexagon.

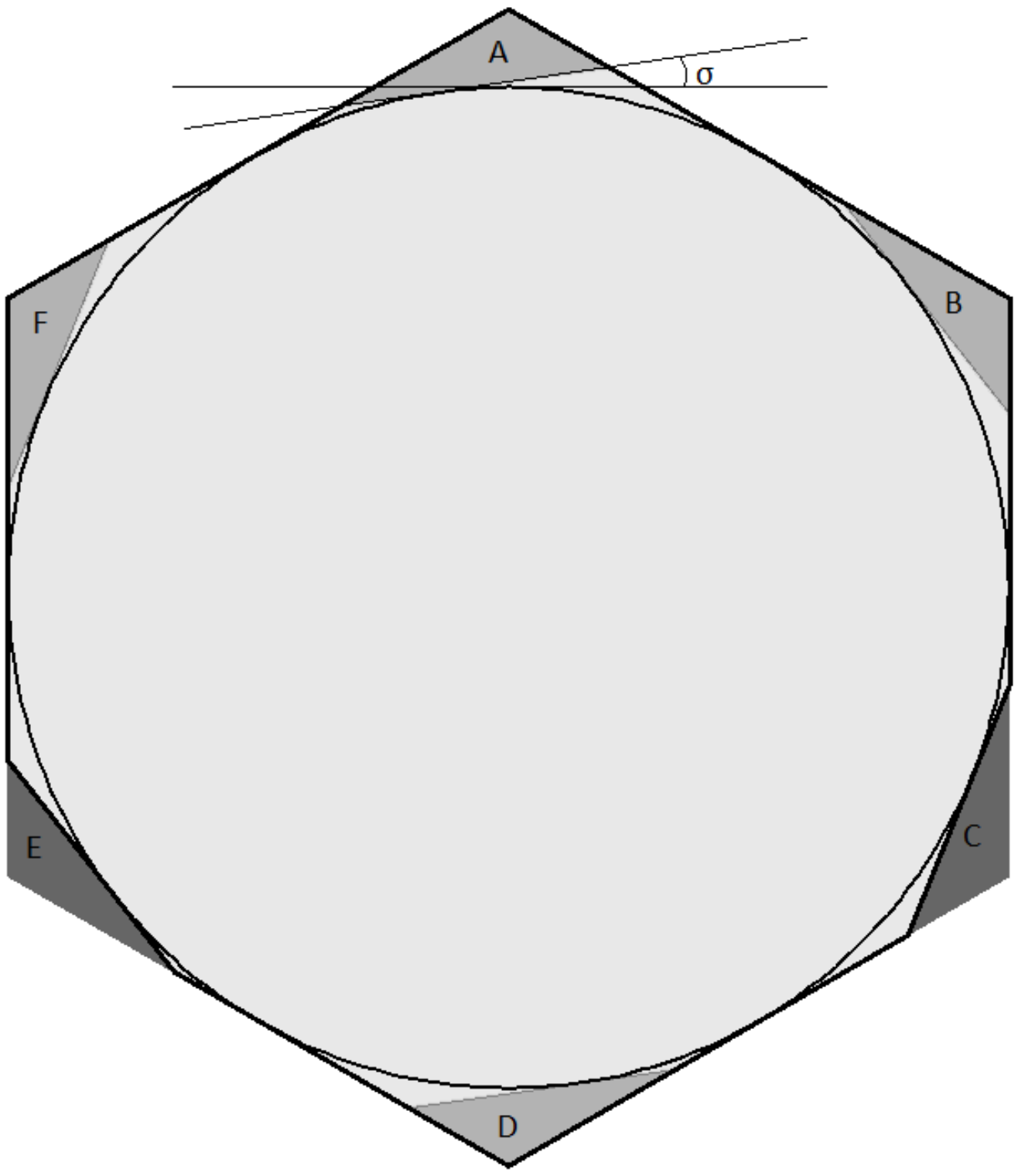

Figure 14: corners cut at slant angle to give a modified Pál cover

Pál's argument can be generalised to show that a universal cover can be formed by cutting off two corners of the regular hexagon tangent to the unit circle but at a different angle to the side of the dodecahedron, provided that this angle is the same for both corners. To show this first consider the regular hexagon with all six corners marked down to a tangent to the circle in such a way that 6-fold rotational symmetry is preserved.

Any shape of diameter one can fit inside the hexagon. Furthermore, it cannot be inside both corner triangles marked A and D since the minimum distance between points in these two regions is one. Similarly it cannot be in both F and C, or both E and B. Making use of the rotational symmetry any combination can be reduced to one of two cases where either the shape is not in corners E, D and C, or not in A, E and C. In either case it is not inside corners E and C so these can be removed leaving a new universal cover.

Once again the minimum area for universal covers of this shape is attained in the case described by Pál where the corners are cut off at the sides of a regular dodecagon, but it is possible that the minimum area for a universal cover within the more general shape is attained for cuts at a different angle. This can be tested using the computational search with a modification whereby any fit inside the hexagon for a shape is rejected if it falls inside the corners that have been sliced off. The results are shown in this graph of area against the angle of the slant away from the dodecahedron edge.

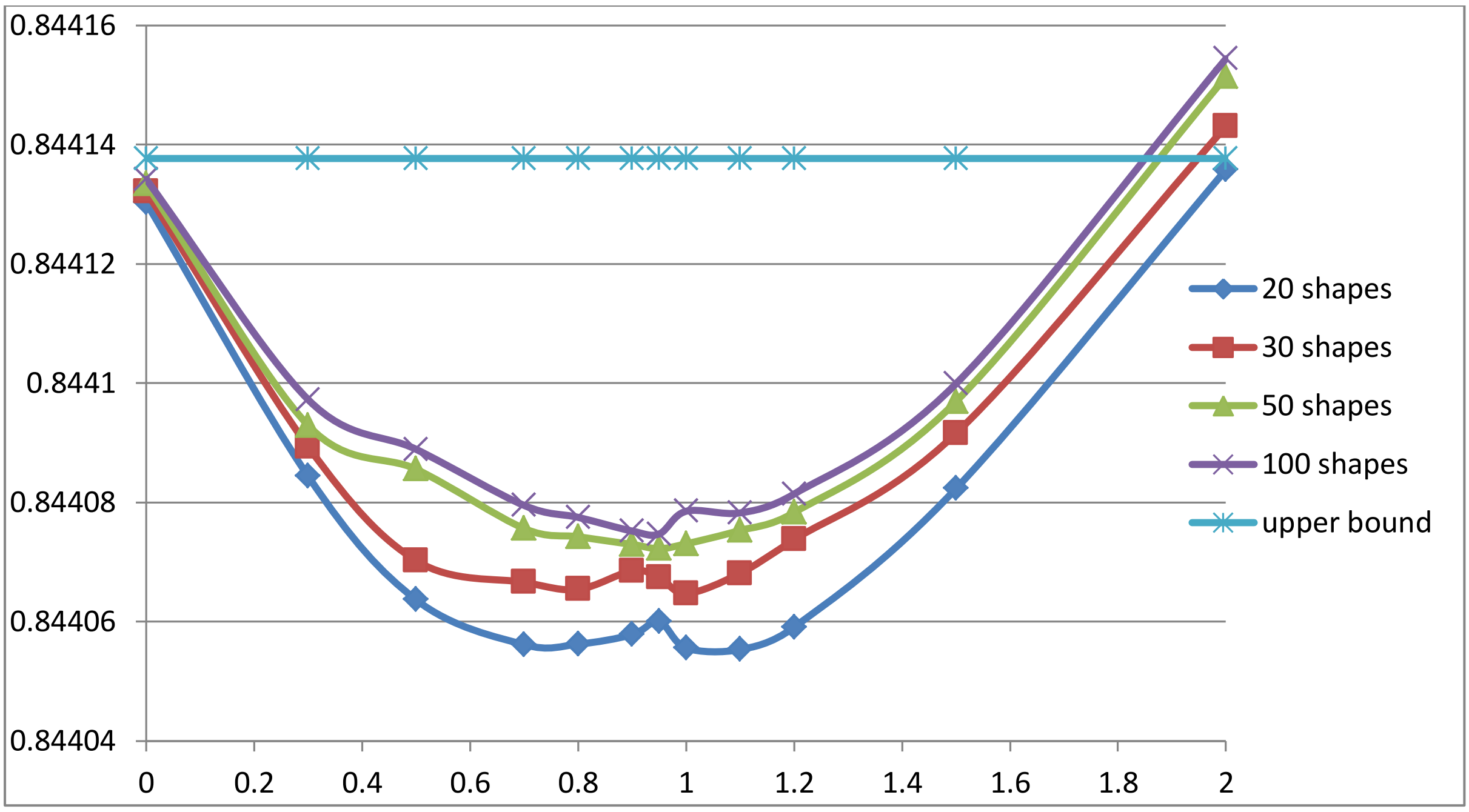

Figure 15: lower bound for modifed Pál hypothesis

This time the results are more encouraging. The minimum is not given by the regular hexagon at slant angle zero. Instead the lowest area is found at around an angle of 1 degree. However, it must be stressed that these points are only lower bounds on the minimum area and the ability to reach the true minimum area under the modified Pál hypothesis could be limited by the choice of possible shapes to include in the search. To be sure of the result it is necessary to prove an upper bound lower than the existing bound provided by Hansen.

Based on the searches for minimum areas inside the regular hexagon a third conjecture emerges

*For a given set of shapes the minimum universal cover inside a regular hexagon is contained within a shape formed by cutting two corners off the hexagon using lines tangent to the circle and at an angle of 120 degrees to each other.*

## 5. Upper Bounds Without Transformations

Any set of points of diameter 1 can be embedded in a curve of constant width one. We also know that any curve of constant width one can fit inside the shape defined by the modified Pál hypothesis. Within that area, there is a smaller set of points defined by the union of all curves of constant width one as they fit inside in any orientation. This shape must be a universal cover but not necessarily the smallest because the minimal cover allows us to choose one of any of the ways the curve fits in. What is the area of the convex hull of this shape as a function of the slant angle?

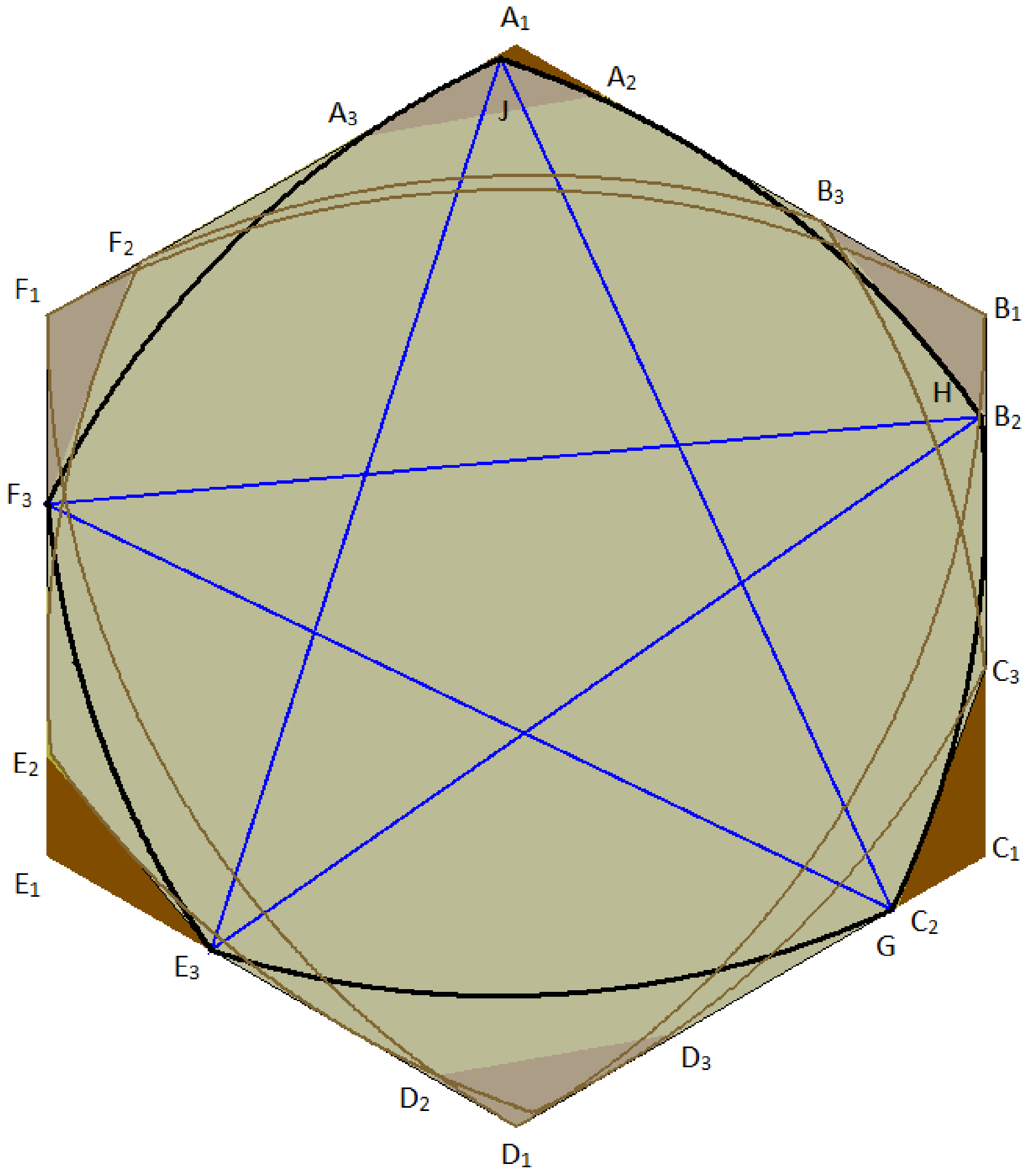

Figure 16: pentagons that define limits of cover

The answer depends on a particular Reuleaux pentagon constructed as follows, using the diagram above and assuming without loss of generality that $A_1A_3$ ia greater than $A_1A_2$. From $F_3$ draw a line of length one to a point G near $C_2$ on the line segment $D_3C_2$. From $F_3$ and $E_3$ draw two lines of length one that meet at a point H near $B_2$. Then from $E_3$ and G draw two lines of length one that meet at a point J near $A_1$ Using the 5 points $F_3$ $E_3$ G, H and J form a Reuleaux pentagon.By construction it fits inside the hexagon and does not enter the regions E, C, F or B. If this pentagon is rotated 180 degrees it forms a second Reuleaux pentagon including $B_3$ and $C_3$ as vertices that fits into the hexagon and also does not enter E, C, F or B. As a third shape add the Reuleaux triangle with vertices at $F_1$, $B_1$ and $D_1$. This shape fits inside the hexagon without entering the triangles at E and C as required.

It can be seen that there are only three small regions that are not covered by this shape. These are near the points $C_2$, $E_2$ and $A_1$ It can be shown that no curve of constant width one within the hexagon that does not enter E and C can enter those three regions.

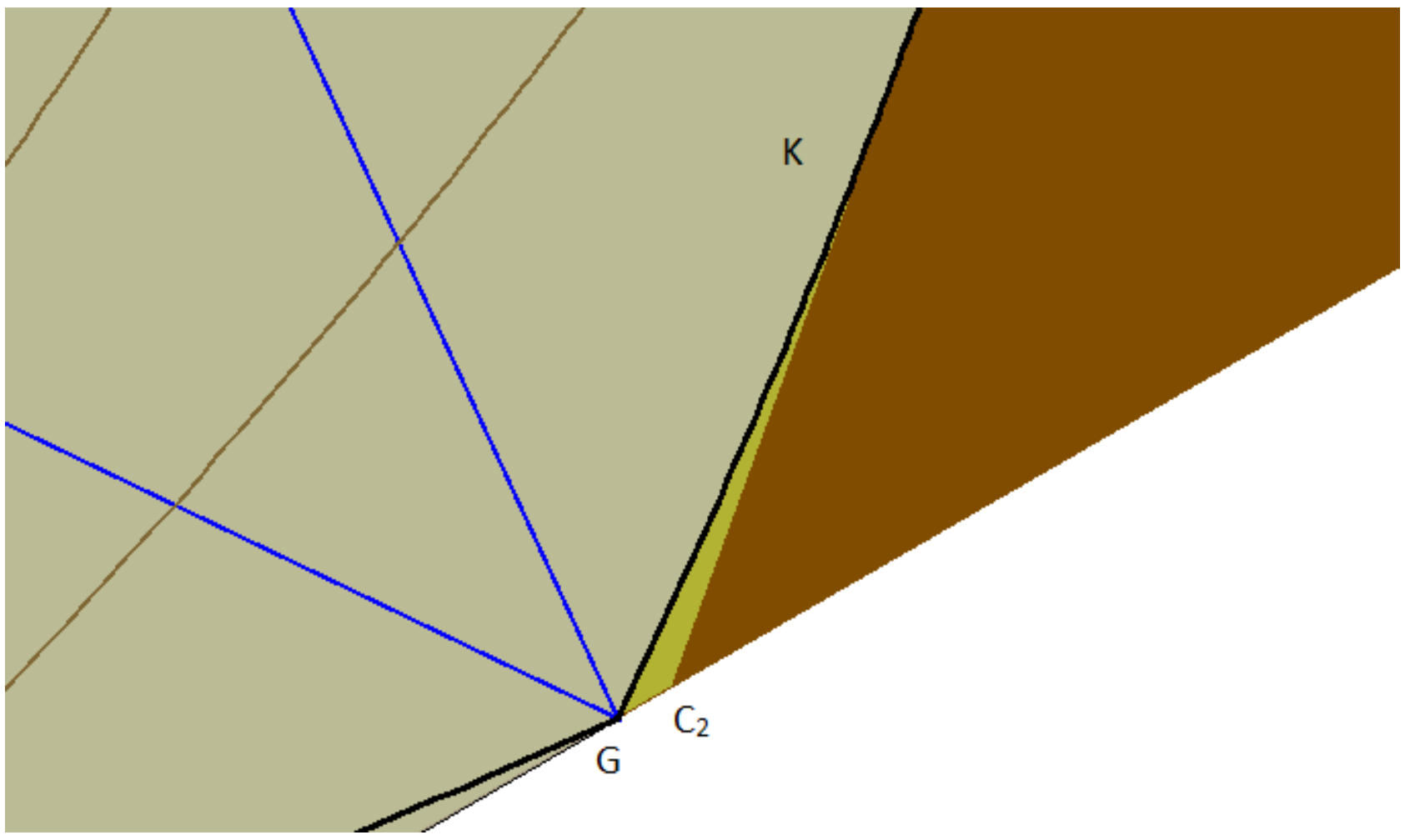

Figure 17: region that can be removed near $C_2$

Firstly the region near $C_2$ is bounded by two straight edges $GC_2$ and $C_2K$ and an arc of radius one KG where K is the point on $C_2C_3$ at distance one from $F_3$. Any curve of unit constant width one fitted into the hexagon without entering triangle C must touch or enter the triangle F. All points in $GC_2K$ are at a distance of one or more from that triangle so no points in the shape can enter this region. In the case of zero slant this region vanishes.

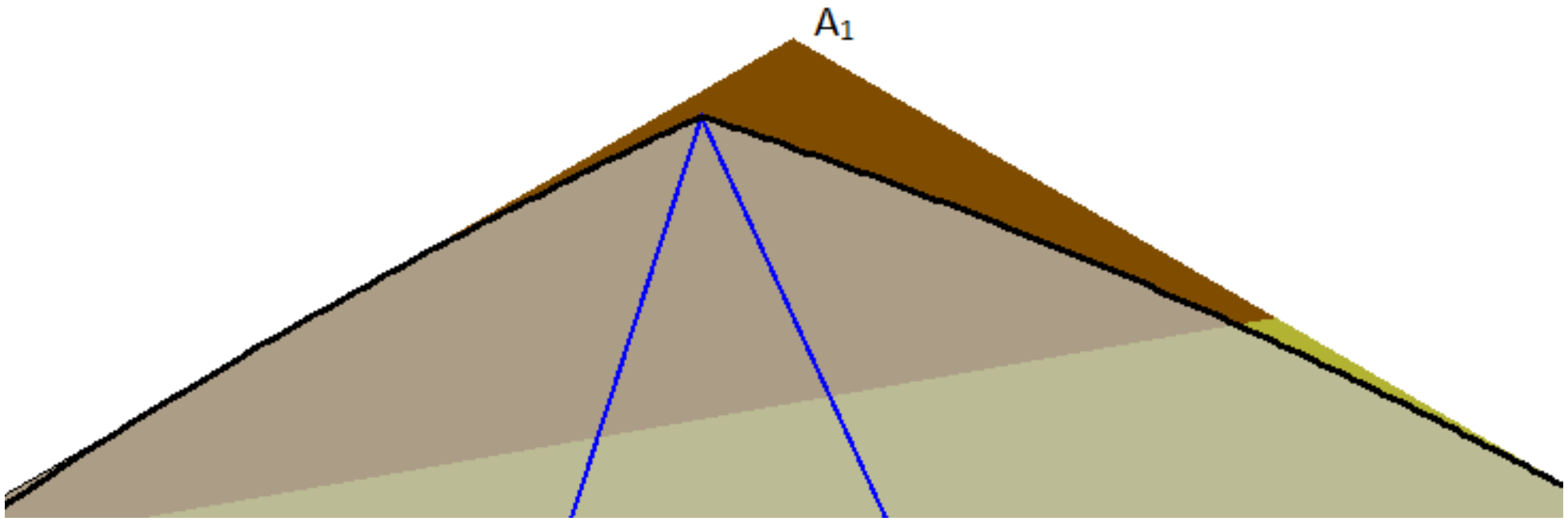

Figure 18: region that can be removed near $A_1$

Secondly the region near $A_1$ is outside the two arcs of radius one centred on $E_3$ and G. A curve of constant with one inside the hexagon must touch each side of the hexagon so it must touch the lines $D_1E_3$ and $D_1G$. This ensures that it cannot enter the region cut off near $A_1$. In the case of zero slant this area is the one removed by Sprague.

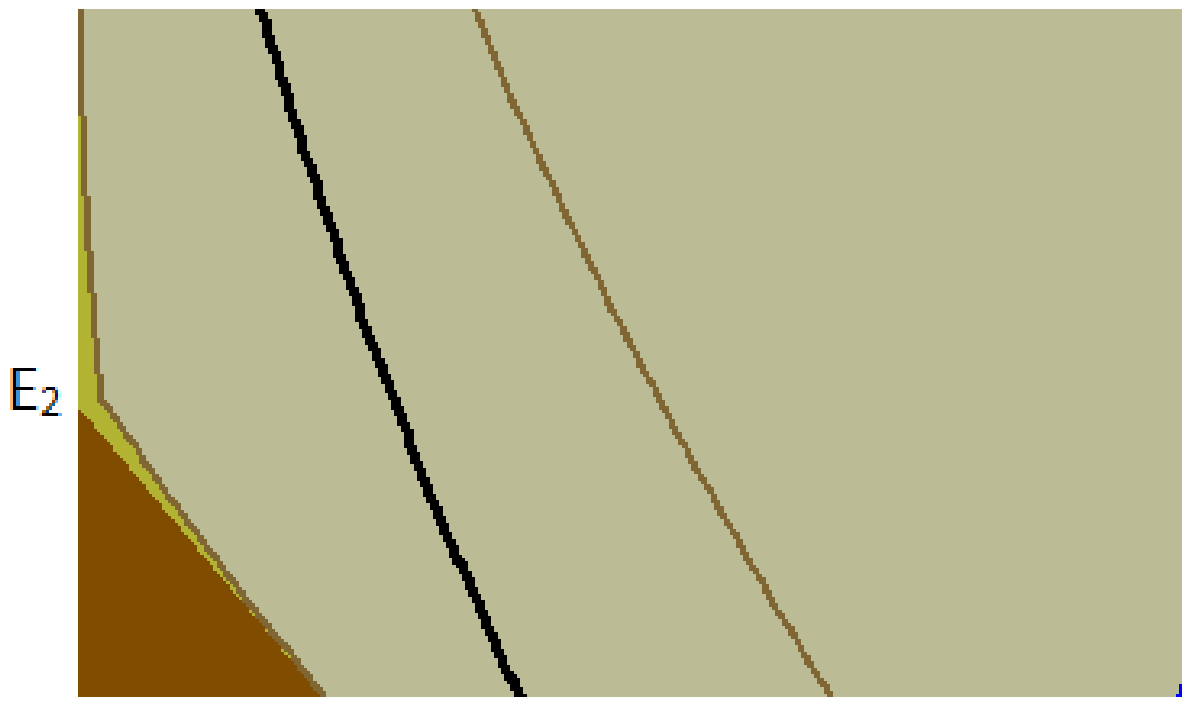

Figure 19: region that can be removed near $E_2$

Thirdly, points in the region cut off around E2 are either greater than distance one from the line segment $B_2C_3$ or greater than one from the triangle B. By similar reasoning there can be no point in this boomerang shaped region that is symmetrical about a line through $E_2$ and $B_2$. This area also vanishes in the case of zero slant $\sigma$.

It now follows that the shape defined by the hexagon with the triangles at E and C and these three smaller regions removed is a universal cover. Its area can be computer and plotted as a function of the angle .

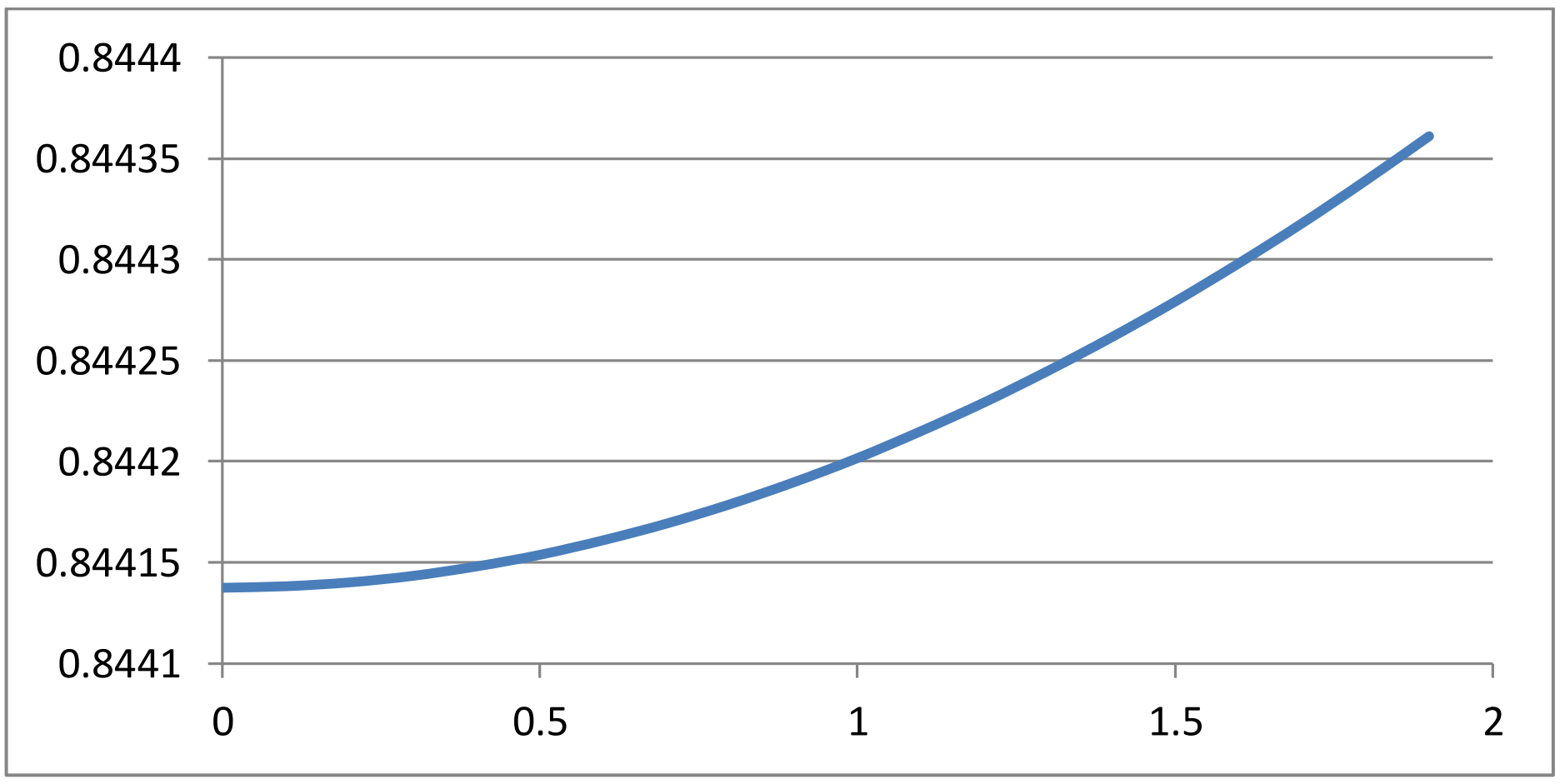

Figure 20: minimum area as a function of slant angle

The minimum of this curve is at $\sigma = 0$ where it reduces to the area determined by Sprague, so this upper bound is not sufficiently good to improve on previous results.

## 6. Upper Bounds Using Rotational Symmetry

To improve the upper bound within a modified Pál hypothesis it is necessary to take into account the ability to rotate or reflect shapes within the area to define a smaller universal cover. Shapes of constant width one fitted within the modified Pál hypothesis can be classified into two types: Those that enter region D and those that don't. If they don't enter region D they will enter A instead. We will call them D-type and A-type. (Those that touch the boundaries of regions A and D are D-type) D-type shapes can be rotated through angles of 120 degree and remain inside the modified Pál hypothesis. This can be used to fit them inside a smaller shape. The A-type shapes can then be treated as before, i.e. construct their convex hull regardless of other positions that may fit. The two shapes defined in this way can be combined to give a new universal cover.

The pentagon $F_3$ $E_3$ G H J is A-type so it is kept in this position. Therefore no improvement in the areas near $E_3$ and $C_2$ are possible using this method. The same pentagon rotated 180 degrees is a D-type fit, so it can be rotated by 120 or 240 degrees. This provides an opportunity to reduce the cover near $E_2$ and $C_3$ where this pentagon was previously a limiting factor. These are the corners where small regions were removed by Hansen in the special case of the regular hexagon so we seek to generalise and improve his results.

An argument that allows small regions near these corners to be removed proceeds as follows. Take a point X near $E_2$. A curve of constant width must touch the line segment $E_3F_2$ at a unique point L. The point at X sets a limit as to how close point L can be to $E_3$ determined by an arc tangent to $E_3F_2$ and passing through X. L cannot be nearer to $E_2$ than the point where this arc touches $E_3F_2$

The shape of width one must also touch $E_3D_2$ at a unique point M. The point X also sets a limit to how close M can be to $E_3$. This uses the observation that if two points X and Y are inside a curve of constant width one then all points on an arc of radius one joining X and Y must also be in the shape. By the modified Pál hypothesis no point on the curve is inside the triangle at E. Therefore if an arc of radius one is drawn from X touching $E_2E_3$ and meeting $E_3D_2$ at M then M can be no closer to $E_2$ than this limit. For some points X such an arc cannot be drawn in which case there is no such restriction on M. Since we are considering only A-type shapes we know that the curve also cannot pass through the triangle D or C so a similar argument can be repeated twice to set limits on the points where the shape touches $D_3C_2$ and then $C_3B_2$ This sets a lower limit to the Point P where the shape touches $C_3B_2$

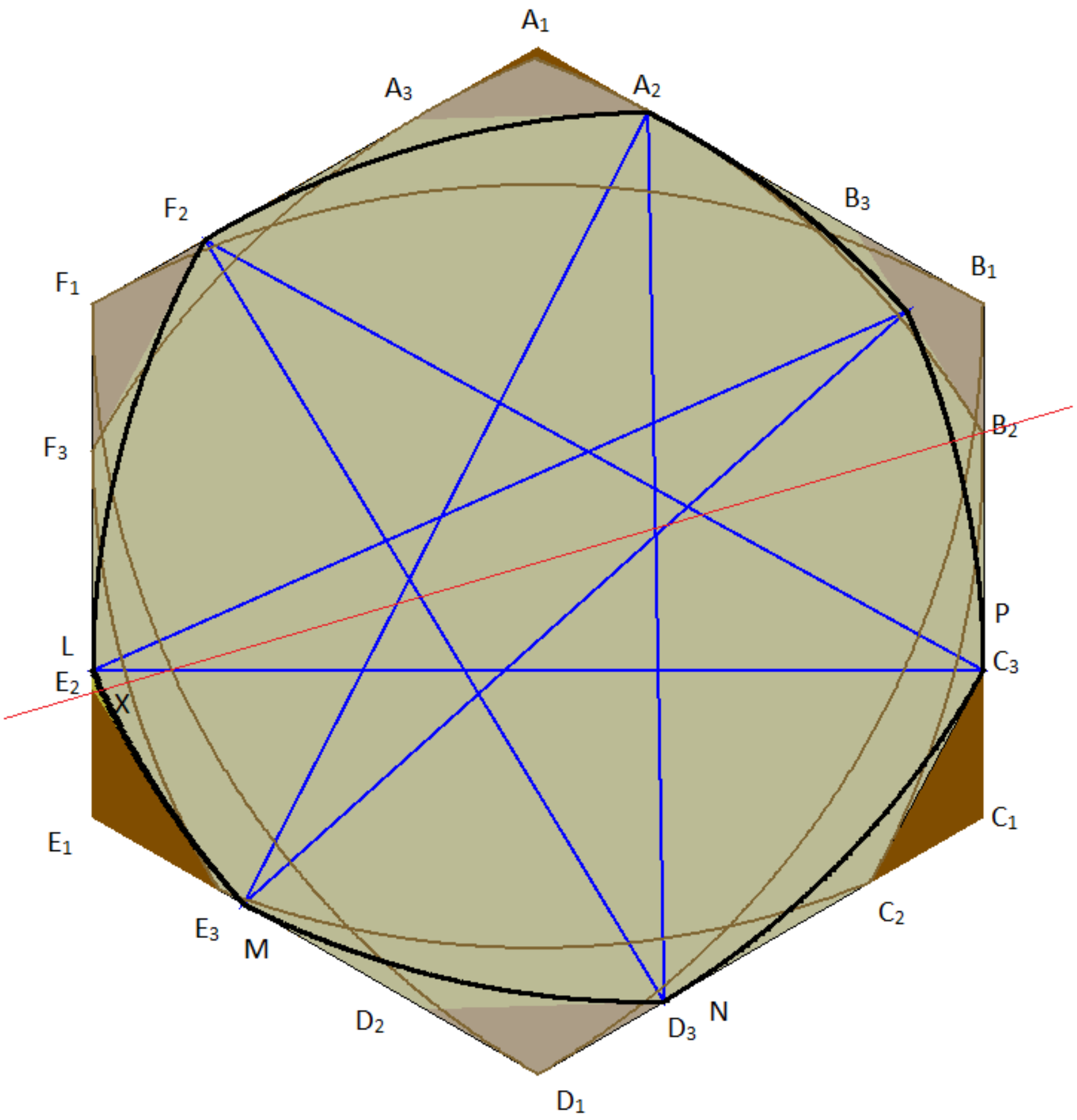

Figure 21: critical heptagons for small slant angles

However the point L where the shape touches F3E2 must be directly opposite the points P where the shape touches $C_3B_2$ The point X sets a lower bound on L and an upper bound on P. If the lower bound is higher than the upper bound there is a contradiction proving that X cannot be inside the an A-type shape. The set of excluded points defines a region that can be removed from the cover of A-type shapes. By a similar argument mirrored a similar but smaller region can be excluded near $C_3$ provided the slant $\sigma$ is not too great.

Taking the case where X is chosen so that it coincides with L and the point P is at the same height leads to a Reuleaux heptagon formed from the arcs and their centres. This defines the upper poiint of the region that can be removed. If this heptagon is reflected about the axis of symmetry through $E_2B_2$ it forms another heptagon which is still inside the Pál hypothesis for small angles. The point where this reflected heptagon touches $E_2E_3$ then defines the lower point of the region.

These two regions do not vanish in the limit of zero slant and are related to areas removed by Hansen in his 1975 paper but are larger. Unfortunately the area is only larger than the previously derived areas above for very small slant angles under 0.0015 degrees. So this cannot account for the minimum lower bounds at around 1 degree.

To complete this analysis the cover formed by D-type fits must also be dealt with. Consider once again the points L, M and N where a D-type curve touches the sides of the hexagon. It will also touch the opposite sides of the hexagon at points L',M' and N' the same distance from the vertices of the hexagon. Taking the three distances $LE_1$=$L'C_1$, $ME_1$=$M'A_1$ and ND1=$N'F_1$ the shape can be rotated through angles of 120 or 240 degrees so that the largest of these three distances is $LE_1$ Since the shape must cross or touch the line $B_3B_2$ and $LE_1 > ME_1$ it follows that a lower bound on the distance $LE_1$ is set when the point L is at distance 1 from the intersection Q of the $B_3B_2$ and the centre line through $E_1B_1$ Furthermore no point near $E_2$ on a D-type shape rotated in this way can be outside the arc of radius one centred on Q. Similarly points near $C_3$ are constrained to be within a distance one of the intersection of $F_3F_2$ and $F_1C_1$. The regions removed from the D-type cover in this way are supersets of the regions removed from the A-type cover. It follows that the region defined by the A-type cover is a full universal cover for all shapes.

## 7. Hansen's Cover and Extension

Before continuing with the investigation of the minimal cover for a modified Pál hypothesis it is worthwhile to revisit the 1992 results of Hansen where small regions near $E_2$ and $C_3$ were removed for the case of zero slant. To accomplish this Hansen used an argument based on the Reuleaux heptagon but with the additional use of reflections. In the case of zero slant the Pál hypothesis is symmetric about the line through $A_1D_1$ so any shape inside the hypothesis can be reflected about that line.

The limiting Heptagon is also symmetric with a vertex on the midpoint of the line $A_3A_2$. Taking the two regions $LE_2R$ and $L'C_3R'$ where R and R' are where this heptagon touches the lines E2E3 and C3C2 respectively, Hansen's argument or a modified version of the arc argument above shows that any A-type shape cannot have points in both of these regions. Therefore the reflections can be used to remove one of these shapes from the cover. On the other side there is still a smaller region that can be removed using the non-symmetric argument as above. Hansen gave a simpler but weaker argument removing just part of this area. Here the area of these regions will be accurately calculated.

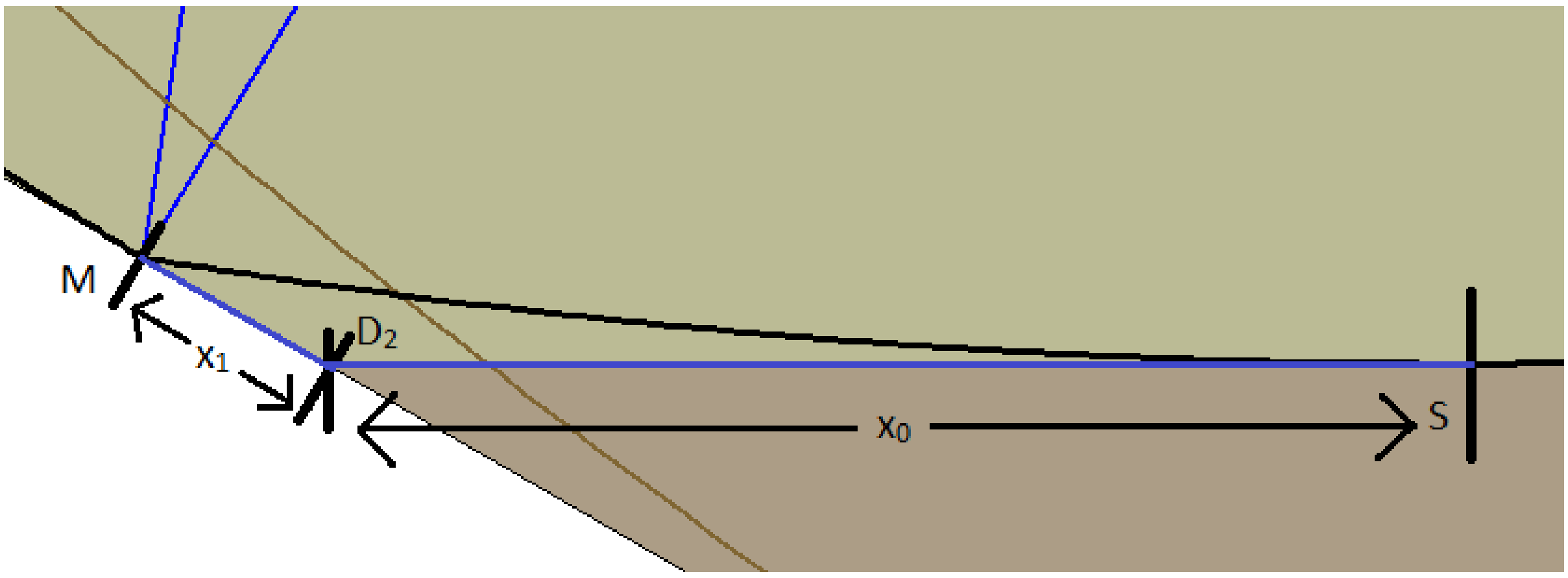

Figure 22: area calculation for regions of regular hexagon

The first step is calculate the position of the vertices M (or N) on the heptagon given that they are equally spaced from the centre line $A_1D_1$ . Let S be the midpoint of $D_2D_3$ and call $x_0$ the distance $D_2S$. From elementary geometry we find $x_0 = 1 - \frac{\sqrt{3}}{2}$ . Knowing $x_0$ the length of $x_1 = MD_2$ is found to be given by

$$x_{i+1} = \frac{2x_i^2}{1 - \sqrt{3}x_i + \sqrt{1 - 2\sqrt{3}x_i - x_i^2}} \qquad (7)$$

A radius of length one from M will meet $B_3A_2$ at the point M' and the line $B_3B_2$ at a new point T. The geometry of the shape $M'B_3T$ is analogous to the geometry of $SD_2M$ but on a smaller scale where the length of $B_3M'$ is given by $x_1$ and the length $x_2 = B_3T$ is given by a second application of the formula above. Once again a radius is drawn from T to intersect the lines meeting at $E_2$. The area defined here is again analogous and the smaller side is $x_3$ calculated using the same formula. This shape is the larger of the two that Hansen removes. To calculate its area use.

$$A_i = \frac{x_i x_{i+1}}{4} - (\theta - sin\theta), \text{ where } \theta = 2sin^{-1}(\frac{d}{2}), \text{ and } d = \sqrt{\frac{x_{i+1}^2}{4} + (x_i + \frac{\sqrt{3}}{2}x_{i+1})^2} \qquad (8)$$

The result is $A_2 = 1.8738 \times 10^{-11}$

The second piece that Hansen removed near $C_3$ is obtained by iterating the construction one more time to give $A_3 = 4.2270 \times 10^{-21}$ .

The area of the boomerang shaped piece that can be removed is $A = 1.3877 \times 10^{-17}$.

## 8. Upper Bound Using Reflections

As shown above, a regular hexagon circumscribed round a circle of diameter one remains a universal cover when two corners are cut off by lines at an angle of 120 degrees to each other and tangent to the circle. These lines can be at a slant angle $\sigma$ to the edges of the regular dodecahedron used by Pál. Three smaller pieces can be cut from this cover near the corners $A_1$, $C_2$ and $E_2$ However, one extra piece near A1 can be removed if we used the freedom to reflect shapes and the remaining convex universal cover then has a minimum area for a non-zero value of the slant angle  . This provides a construction for a convex universal cover that has a smaller area than previous best results.

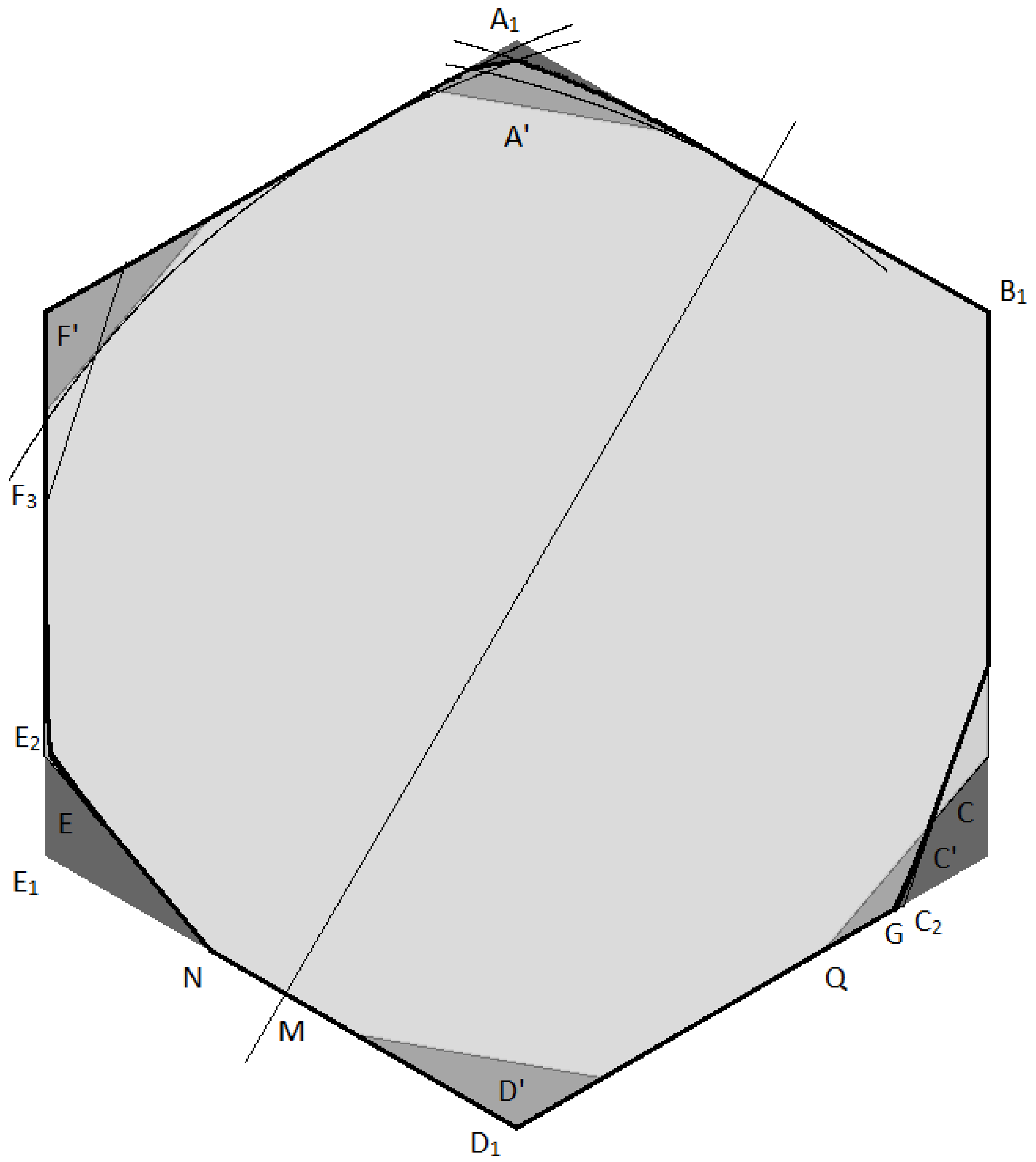

Figure 23: construction of universal cover using reflections

The axis of reflection to use is the line from the midpoint M of the side of the hexagon from $E_1$ to $D_1$, to the midpoint of the opposite side. A shape fitted into the hexagon with the corners E and C removed can be reflected about this axis provided it also does not enter the triangles F' and D' which are the reflections of C and E about the axis. When this is the case we will choose to reflect the shape if it touches the side from E to D at a point nearer to $E_1$ than $D_1$. Remember that it touches the opposite side at the opposite point which is therefore also reflected to be nearer the corner of the hexagon at B than the one at A. If we draw an arc centred on M of radius one it cuts off a larger region near A and no shape that can be reflected can be beyond this line. The point where this meets the arc centred on G is marked W. Recall that G is the point at distance one from $F_3$, so no point inside any covered shape can lie outside the arc centred on G.

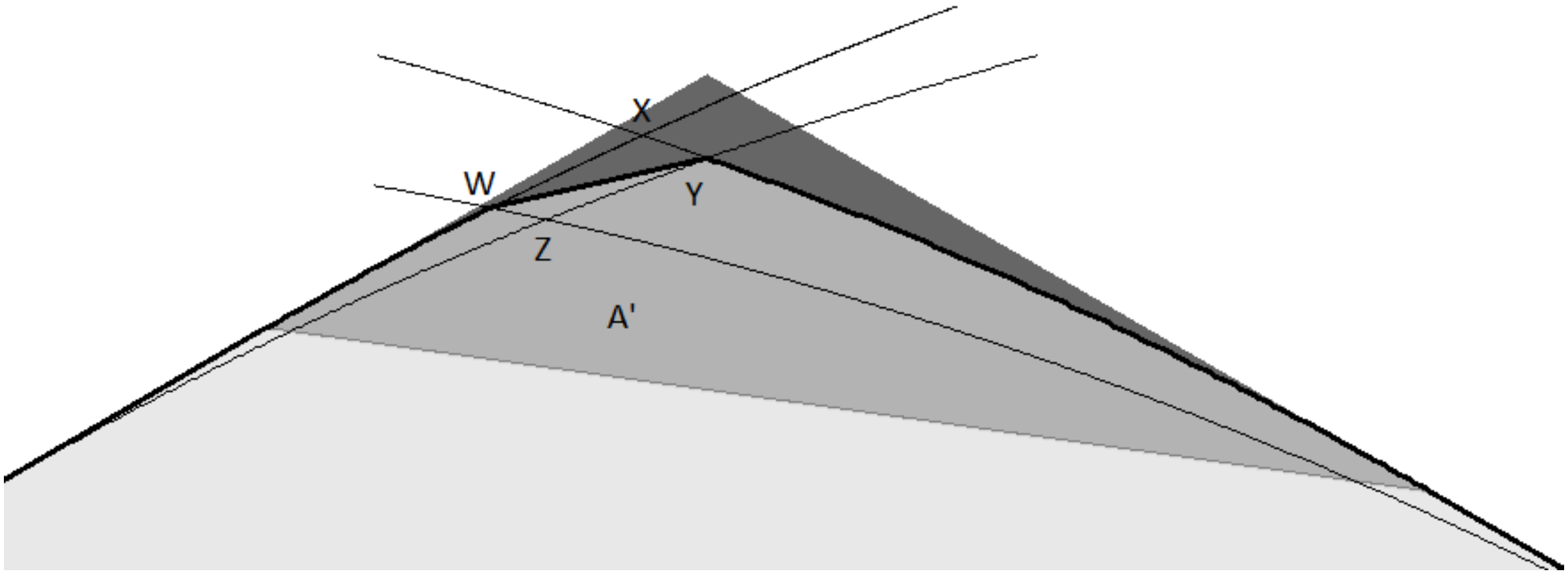

Figure 24: region XYZW that can be removed

The only shapes that can have points outside that arc centred on M are ones that cannot be reflected. Therefore they must enter some point inside the regions F' or D' If they have a point in D' then they cannot have a point in A' which is the triangle whose points are at a distance of more than one from all points in D'. Draw one more arc centred on Q at the corner of the triangle C' which is the reflection of the region F. All points in C' are at a distance of one or more than one from points in F' so the arc will touch the region F' but not enter it. This arc will meet the arc centred on M at a point Z and the arc centred on N at a point Y. Now consider the fate of points inside the region XYZW bounded by the four arcs. It is a general property of curves of constant width one that if two points are inside the curve then all points on an arc of radius one through the two points are also inside the curve. Suppose then that a shape fitted inside the hexagon had a point in XYZW and also in F' We could then join those two points with an arc but between the two points it would be outside the arc centred on Q and would therefore go outside the hexagon. This is in contradiction with the premise so we conclude that no shape fitted in the hexagon can have a point in both XYZW and F'. It can also be verified that for angles $\sigma$ less than 9 degrees the region XYZW is inside the triangle A'. Therefore shapes with a point inside XYZW do not have points in F' or D' and can be reflected. However, we have already determined that such shapes will not have points in this region. This proves that the region XYZW can be removed from the universal cover.

It turns out that this is now sufficient to construct a universal cover smaller than the ones of Hansen and Duff. Even if we restrict ourselves to the convex case and remove only the part of this region that leaves a convex shape, the area of the universal cover for an angle $\sigma = 0.52$ degrees can be computed to be 0.8441121.

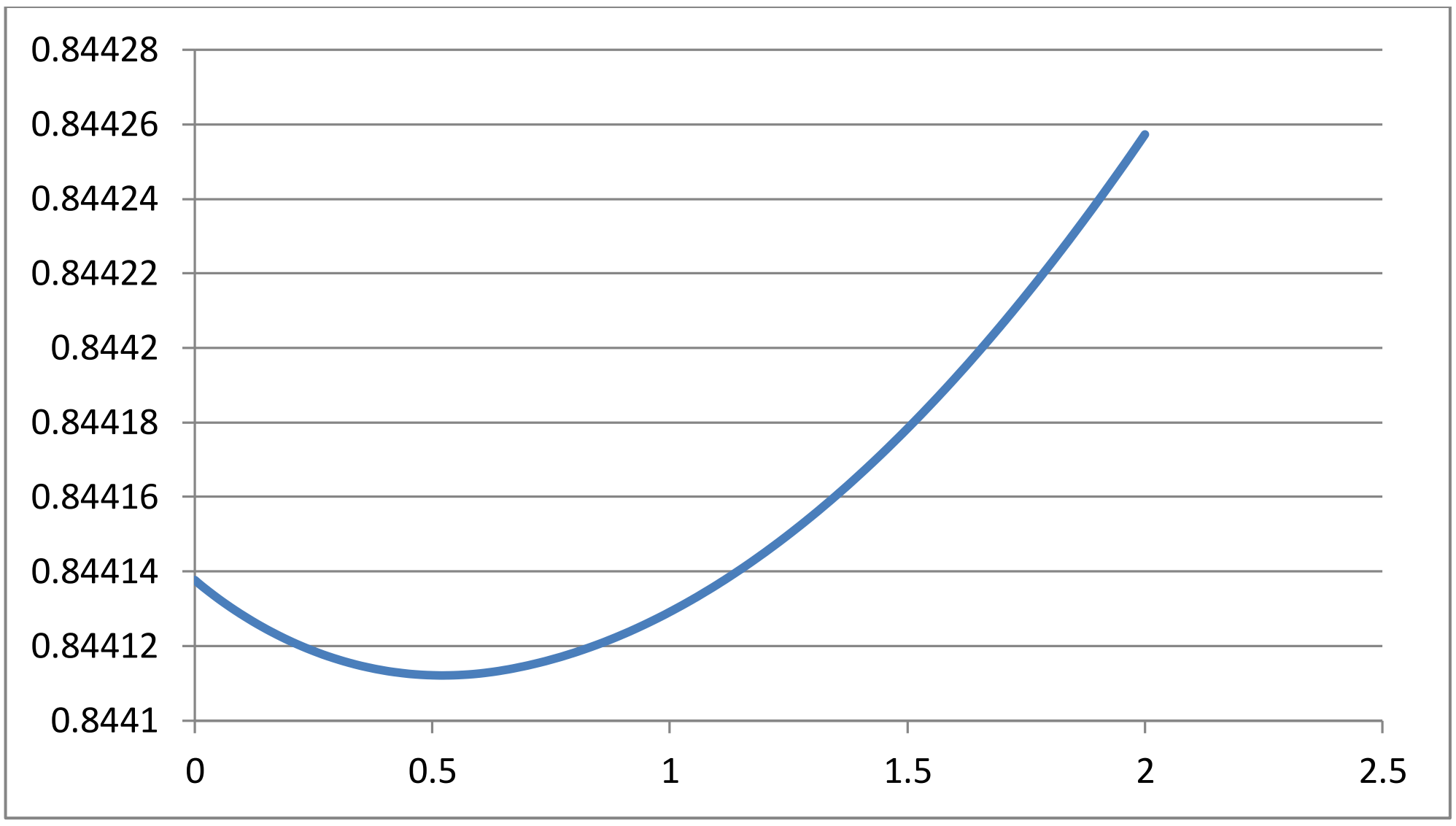

Figure 25: area as a function of slant angle using reflections

## 9. Opportunity for further progress

The upper bound of 0.8441121 at a slant angle of 0.54 degrees is not the final answer. The computational searches using a regular hexagon and a fixed slant angle gave a minimum area of about 0.84408 at a slant angle of about 1.0 degrees. If we accept the conjectures based on the simulations then the minimum cover should be enclosed in a regular hexagon with a modified Pál hypothesis. This result is therefore likely to be closer to the final answer, but where was more area removed?

To answer this, the simulations were ran again for a slant angle of 1 degree with background masks showing the regions removed in the reflection analysis. This shows that there are too areas where more area can be removed.

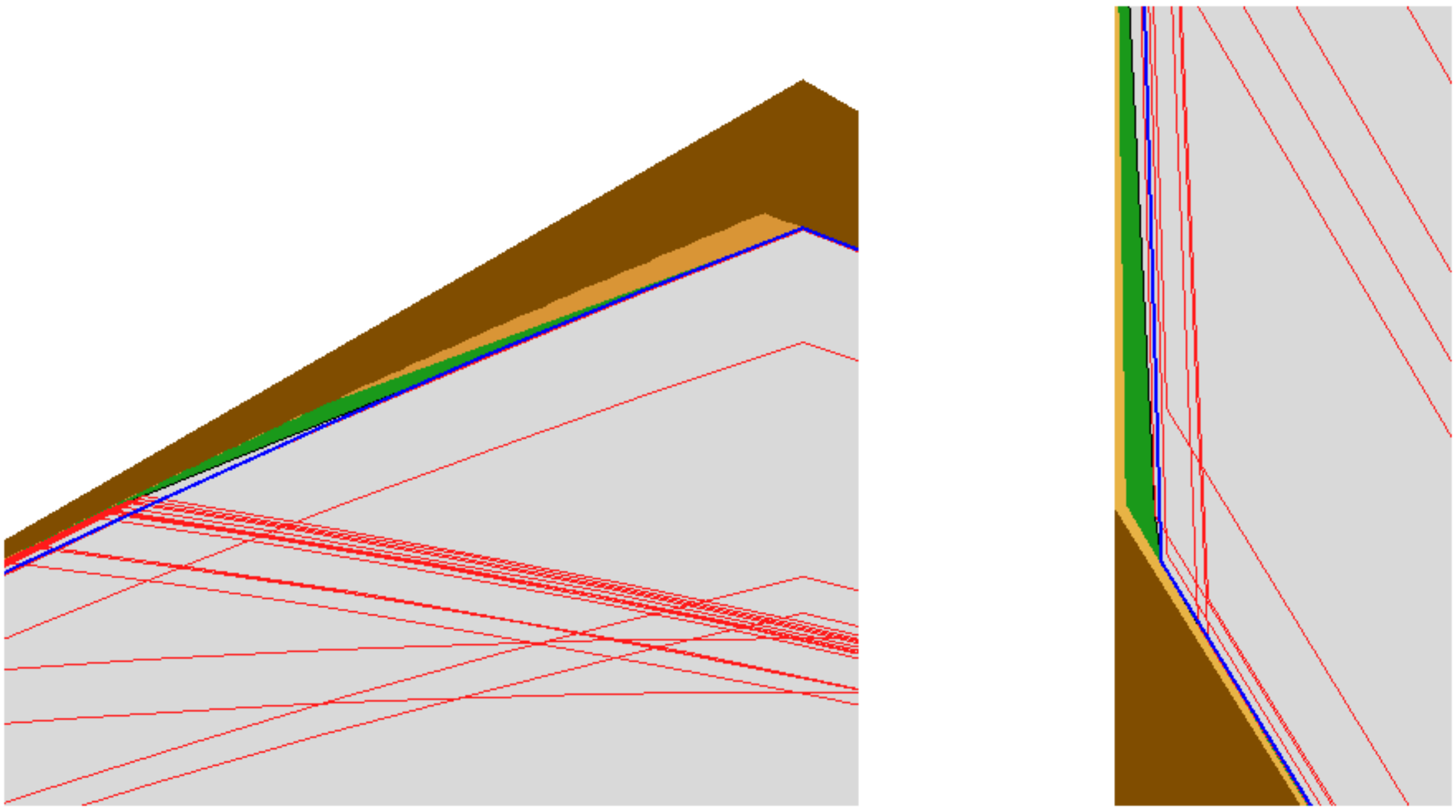

Figure 26: orange is area already removed, green is further area that can be removed

## 10. Summary and Conclusions

Computational methods have been used to investigate Lebesgue's Universal Covering Problem. This gave indications of what the minimal cover is like and once known how to proceed, it was then easy to prove an upper bound for the convex cover improving previous results.

The conclusions can be framed in a series of conjectures. First it was noticed from simulated annealing that for any set of shapes of constant width one the minimal convex covering can be contained within a hexagon with opposite sides parallel and distance one apart. It would be a reasonable conjecture that this is indeed the case.

If the set of shapes includes the circle then we know more specifically that the hexagon can be circumscribed round a circle with parallel opposite sides. Any such hexagon can be shown itself to be a universal covering so there is a minimal convex covering that it contains. The minimal such cover for a set of shapes can be computed efficiently and the results show that the cover with the minimal area cannot differ by much from the case of the regular hexagon. Again it would be a reasonable conjecture that the regular hexagon is precisely the minimal case for all such hexagons.

Given a regular hexagon it was observed that the minimal convex cover within the hexagon for any set of shapes of constant width can be contained in a shape formed by cutting off two of the corners using two lines tangent to the inscribed circle at 120 degrees to each other (a modified Pál hypothesis) Again it could be conjectured that this is the case. Computations of the minimal area within a modified Pál hypothesis indicated that the minimum convex area is given when the lines are at a small angle away from the edges of the regular dodecahedron. It was then easy to prove that a convex universal covering could be constructed in this way improving previous upper bounds.

Further work would probably reduce the area still further and may lead to a conjectured minimum convex cover. A strategy to prove this would be to prove the three individual conjectures and then show that the conjectured minimal convex cover is the minimum within a modified Pál hypothesis. None of these steps looks simple.

The case of non-convex covers has not been investigated in the same way but it is easy to see that a further reduction in the area can be attained for the non-convex case.

## Acknowledgments

My thanks go to John Baez for bringing attention to this problem via his blog Azimuth, and for providing encouragement and comments on this paper.

## Note on Publication

It would have been useful to subject this paper to peer-review so I have searched for suitable journals that meet the following conditions:

- They are Open Access,
- free of author charges or allow waivers for unfunded research,
- do not require a TeX rewrite and
- accept papers from authors without an academic affiliation.

Unfortunately I could not find any journal of geometry that satisfies all of these needs. I therefore apologise for any faults this article may contain and would be grateful for error reports from anyone who finds them.

## References

[1] H Lebesgue, Letter to J. Pál (1914) as quoted in [3].

[2] Paul J. Kelly; Max L. Weiss (1979). Geometry and Convexity: A Study in Mathematical Methods. Wiley. pp. Section 6.4.

[3] Pal, J., 'Über ein elementares Variationsproblem', Danske Mat.-Fys. Meddelelser III, 2 (1920).

[4] Sprague, R., 'Über ein elementares variationsproblem', Mat. Tidsskrift (1936), 96-99.

[5] Hansen, H. C., 'A small universal cover of figures of unit diameter', Geom. Dedicata 4 (1975), 165-172.

[6] Hansen, H. C., "Towards the minimal universal cover", Normat 29 (1981) 115-119, 148

[7] Hansen, H. C., 'Small universal covers for sets of unit diameter', Geom. Dedicata 42 (1992) 205-213

[8] Duff, G. F. D., 'A smaller universal cover for sets of unit diameter', C. R. Math. Acad. Sci. 2

(1980), 37-42

[9] Rennie B.C. "The Search for a Universal Cover", Eueka 3 (Ottawa) No 3, (March 1977)

[10] Duff. G.F.D., "On Universal Covering Sets and Translation Covers in the Plane", James Cook Mathematical Notes, 23, Vol 2 (1980) p109-201

[11] Elekes, G. "Generalized Breadths, Circular Cantor Sets, and the Least Area UCC", Discrete & Computational Geometry (1994), Volume 12, Issue 1, pp 439-449

[12] P. Brass, M. Sharifi: "A lower bound for Lebesgue's universal cover problem", International Journal on Computational Geometry & Applications 15 (2005) 537—544

[13] P. Brass, W. O. J. Moser, J. Pach: "Research problems in discrete geometry." Springer 2005,

[14] B. Grunbaum, "Borsuk's Problem and Related Questions", in Proceedings Symposia in Pure mathematics, Vol VII, Convexity, AMS, Providence (1963) 271-284